\documentclass[11pt]{article}
\usepackage{amsmath,amssymb,amsthm,geometry,verbatim,enumerate,float,tikz}
\usetikzlibrary{decorations,arrows}
\usepgflibrary{decorations.pathreplacing} 
\newtheorem{thm}{Theorem}
\newtheorem{lemma}[thm]{Lemma}
\newtheorem{coro}[thm]{Corollary}

\newtheorem{claim}{Claim}

\usepackage{setspace}

\begin{document}

\onehalfspace

\title{A Characterization of Mixed Unit Interval Graphs}

\author{Felix Joos}

\date{}

\maketitle

\begin{center}
Institut f\"{u}r Optimierung und Operations Research, 
Universit\"{a}t Ulm, Ulm, Germany\\
\texttt{felix.joos@uni-ulm.de}
\end{center}

\begin{abstract}
We give a complete characterization of mixed unit interval graphs,
the intersection graphs of closed, open, and half-open unit intervals of the real line.
This is a proper superclass of the well known unit interval graphs.
Our result solves a problem posed by Dourado, Le, Protti, Rautenbach and Szwarcfiter  
(Mixed unit interval graphs, Discrete Math. \textbf{312}, 3357-3363 (2012)).

\bigskip

\noindent {\bf Keywords:} unit interval graph; proper interval graph; intersection graph
\end{abstract}

\section{Introduction}
A graph $G$ is an \textit{interval graph}, if there is a function $I$ from the vertex set of $G$
to the set of intervals of the real line such that two vertices are adjacent if and only if
their assigned intervals intersect.
The function $I$ is an \textit{interval representation} of $G$.
Interval graphs are well known and investigated \cite{fishburn, golumbic, LekBo}.
There are several different algorithms that decide, if a given graph is an interval graph. 
See for example \cite{corneiletal}.

An important subclass of interval graphs are unit interval graphs.
An interval graph $G$ is a \textit{unit interval graph},
if there is an interval representation $I$ of $G$ such that
$I$ assigns to every vertex a closed interval of unit length.
This subclass is well understood and easy to characterize structurally \cite{roberts}
as well as algorithmically \cite{corneil}.

Frankl and Maehara \cite{fm} showed that it does not matter,
if we assign the vertices of $G$ only to closed intervals or only to open intervals of unit length.
Rautenbach and Szwarcfiter \cite{rs} characterized, by a finite list of forbidden induced subgraphs,
all interval graphs $G$ such that there is an interval representation of $G$ that uses only open and closed unit intervals.

Dourado et al. \cite{dlprs} gave a characterization of all diamond-free interval graphs that have an interval representation 
such that all vertices are assigned to unit intervals,
where all kinds of unit intervals are allowed and a diamond is a complete graph on four vertices minus an edge.
Furthermore, they made a conjecture concerning the general case.
We prove that their conjecture is not completely correct and give a complete characterization of this class.
Since the conjecture is rather technical and not given by a list of forbidden subgraphs, we refer the reader to \cite{dlprs} for a detailed formulation of the conjecture,
but roughly speaking, they missed the class of forbidden subgraphs shown in Figure \ref{graphsT}.

In Section 2 we introduce all definitions and relate our result to other work.
In Section 3 we state and prove our results.

\section{Preliminary Remarks}
We only consider finite, undirected, and simple graphs.
Let $G$ be a graph.
We denote by $V(G)$ and $E(G)$ the vertex and edge set of $G$, respectively.
If $C$ is a set of vertices,
then we denote by $G[C]$ the subgraph of $G$ induced by $C$.
Let $\mathcal{M}$ be a set of graphs.
We say $G$ is $\mathcal{M}$-\textit{free},
if for every $H\in\mathcal{M}$, the graph $H$ is not an induced subgraph of $G$.
For a vertex $v\in V(G)$,
let the \textit{neighborhood} $N_G(v)$ of $v$ be the set of all vertices that are adjacent to $v$ and
let the \textit{closed neighborhood} $N_G[v]$ be defined by $N_G(v)\cup \{v\}$.
Two distinct vertices $u$ and $v$ are \textit{twins} (in $G$)
if $N_G[u]=N_G[v]$.
If $G$ contains no twins, then $G$ is \textit{twin-free}.

Let $\mathcal{N}$ be a family of sets.
We say a graph $G$ has an $\mathcal{N}$-\textit{intersection representation},
if there is a function $f:V(G)\rightarrow \mathcal{N}$
such that for any two distinct vertices $u$ and $v$
there is an edge joining $u$ and $v$ if and only if $f(u)\cap f(v)\not=\emptyset$.
If there is an $\mathcal{N}$-intersection representation for $G$,
then $G$ is an $\mathcal{N}$-\textit{graph}.
Let $x,y\in \mathbb{R}$.
We denote by $$[x,y]=\{z\in \mathbb{R}: x\leq z\leq y\}$$ the \textit{closed interval},
by $$(x,y)=\{z\in \mathbb{R}: x< z<y\}$$ the \textit{open interval},
by $$(x,y]=\{z\in \mathbb{R}: x< z\leq y\}$$ the \textit{open-closed interval}, and
by $$[x,y)=\{z\in \mathbb{R}: x\leq z< y\}$$ the \textit{closed-open interval} of $x$ and $y$.
For an interval $A$, 
let $\ell(A)=\inf\{x\in \mathbb{R}:x\in A\}$ and $r(A)=\sup\{x\in \mathbb{R}:x\in A\}$.
If $I$ is an interval representation of $G$ and $v\in V(G)$,
then we write $\ell(v)$ and $r(v)$ instead of $\ell(I(v))$ and $r(I(v))$, respectively,
if there are no ambiguities.
Let
$\mathcal{I}^{++}$ be the set of all closed intervals,
$\mathcal{I}^{--}$ be the set of all open intervals,
$\mathcal{I}^{-+}$ be the set of all open-closed intervals,
$\mathcal{I}^{+-}$ be the set of all closed-open intervals, and
$\mathcal{I}$ be the set of all intervals.
In addition, let
$\mathcal{U}^{++}$ be the set of all closed unit intervals,
$\mathcal{U}^{--}$ be the set of all open unit intervals,
$\mathcal{U}^{-+}$ be the set of all open-closed unit intervals,
$\mathcal{U}^{+-}$ be the set of all closed-open unit intervals, and
$\mathcal{U}$ be the set of all unit intervals.
We call a $\mathcal{U}$-graph a \textit{mixed unit interval graph}.

By a result of \cite{dlprs} and \cite{rs}, every interval graph is an $\mathcal{I}^{++}$-graph.
With our notation unit interval graphs equals $\mathcal{U}^{++}$-graphs.
An interval graph $G$ is a \textit{proper interval graph}
if there is an interval representation of $G$ such that $I(u)\not\subseteq I(v)$ for every distinct $u,v\in V(G)$.

The next result due to Roberts characterizes unit interval graphs.

\begin{thm}[Roberts \cite{roberts}] \label{thmroberts}
The classes of unit interval graphs, proper interval graphs, and $K_{1,3}$-free interval graphs are the same.
\end{thm}

The second result shows that several natural subclasses of mixed unit interval graphs actually coincide with the class of unit interval graphs.

\begin{thm}[Dourado et al., Frankl and Maehara \cite{dlprs,fm}]
The classes of $\mathcal{U}^{++}$-graphs,
$\mathcal{U}^{--}$-graphs,
$\mathcal{U}^{+-}$-graphs,
$\mathcal{U}^{-+}$-graphs, and
$\mathcal{U}^{+-}\cup\mathcal{U}^{-+}$-graphs are the same.
\end{thm}

A graph $G$ is a \textit{mixed proper interval graph} (respectively an \textit{almost proper interval graph})
if $G$ has an interval representation 
$I:V(G)\rightarrow \mathcal{I}$  (respectively $I:V(G)\rightarrow \mathcal{I}^{++}\cup \mathcal{I}^{--}$)
such that
\begin{itemize}
	\item there are no two distinct vertices $u$ and $v$ of $G$
	with $I(u),I(v)\in \mathcal{I}^{++}$,
	$I(u)\subseteq I(v)$, and $I(u)\not=I(v)$, and
	\item for every vertex $u$ of $G$ with $I(u)\notin \mathcal{I}^{++}$,
	there is a vertex $v$ of $G$ with $I(v)\in \mathcal{I}^{++}$,
	$\ell(u)=\ell(v)$, and $r(u)=r(v)$.
\end{itemize}

A natural class extending the class of unit interval graphs are $\mathcal{U}^{++}\cup\mathcal{U}^{--}$-graphs.
These were characterized by Rautenbach and Szwarcfiter.

\begin{figure}[t]
\begin{center}
\begin{tikzpicture}[scale=1]
\def\ver{0.1} 
\def\x{1}

\def\xa{0}
\def\ya{0}

\def\xb{4}
\def\yb{0}

\def\xc{8}
\def\yc{0}

\def\xd{12}
\def\yd{0}

\path[fill] (\xa,\ya) circle (\ver);
\path[fill] (\xa+1,\ya) circle (\ver);
\path[fill] (\xa+2,\ya) circle (\ver);
\path[fill] (\xa+3,\ya) circle (\ver);
\path[fill] (\xa+1.5,\ya+1.5) circle (\ver);

\draw[thick] (\xa,\ya)--(\xa+1.5,\ya+1.5)
(\xa+1,\ya)--(\xa+1.5,\ya+1.5)
(\xa+2,\ya)--(\xa+1.5,\ya+1.5)
(\xa+3,\ya)--(\xa+1.5,\ya+1.5);

\node (1) at (\xa+1.5,\ya-0.7) {$K_{1,4}$};

\path[fill] (\xb,\yb) circle (\ver);
\path[fill] (\xb+1,\yb) circle (\ver);
\path[fill] (\xb+2,\yb) circle (\ver);
\path[fill] (\xb+3,\yb) circle (\ver);
\path[fill] (\xb+1.5,\yb+1.5) circle (\ver);

\draw[thick] (\xb,\yb)--(\xb+1.5,\yb+1.5)
(\xb+2,\yb)--(\xb+1,\yb)--(\xb+1.5,\yb+1.5)
(\xb+2,\yb)--(\xb+1.5,\yb+1.5)
(\xb+3,\yb)--(\xb+1.5,\yb+1.5);

\node (1) at (\xb+1.5,\yb-0.7) {$K_{1,4}^*$};

\path[fill] (\xc,\yc) circle (\ver);
\path[fill] (\xc+1.5,\yc) circle (\ver);
\path[fill] (\xc+3,\yc) circle (\ver);
\path[fill] (\xc+1,\yc+1.5) circle (\ver);
\path[fill] (\xc+2,\yc+1.5) circle (\ver);

\draw[thick] (\xc,\yc)--(\xc+1,\yc+1.5)--(\xc+1.5,\yc)--(\xc+2,\yc+1.5)--
(\xc+3,\yc)--(\xc+1,\yc+1.5)--(\xc+2,\yc+1.5)--(\xc,\yc);

\node (1) at (\xc+1.5,\yc-0.7) {$K_{2,3}^*$};

\path[fill] (\xd,\yd) circle (\ver);
\path[fill] (\xd+1,\yd) circle (\ver);
\path[fill] (\xd+2,\yd) circle (\ver);
\path[fill] (\xd+3,\yd) circle (\ver);
\path[fill] (\xd+1,\yd+1.5) circle (\ver);
\path[fill] (\xd+2,\yd+1.5) circle (\ver);

\draw[thick] (\xd,\yd)--(\xd+1,\yd+1.5)--
(\xd+1,\yd)--(\xd+2,\yd+1.5)--
(\xd+2,\yd)--(\xd+1,\yd+1.5)--
(\xd+2,\yd+1.5)--(\xd+3,\yd);

\node (1) at (\xd+1.5,\yd-0.7) {$K_{2,4}^*$};

\end{tikzpicture}
\end{center}
\caption{Forbidden induced subgraphs for twin-free $\mathcal{U}^{++}\cup\mathcal{U}^{--}$-graphs.}\label{u+-graph}
\end{figure}
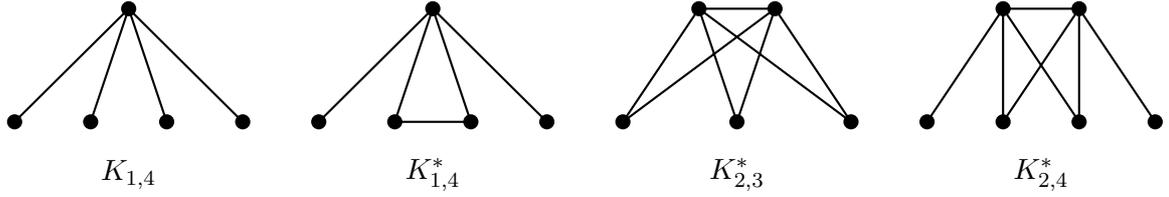

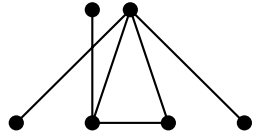
\begin{figure}[t]
\begin{center}
\begin{tikzpicture}[scale=1]
\def\ver{0.1} 
\def\x{1}

\def\xa{0}
\def\ya{0}

\def\xb{4}
\def\yb{0}

\def\xc{8}
\def\yc{0}

\def\xd{12}
\def\yd{0}

\path[fill] (\xb,\yb) circle (\ver);
\path[fill] (\xb+1,\yb) circle (\ver);
\path[fill] (\xb+2,\yb) circle (\ver);
\path[fill] (\xb+3,\yb) circle (\ver);
\path[fill] (\xb+1.5,\yb+1.5) circle (\ver);
\path[fill] (\xb+1,\yb+1.5) circle (\ver);

\draw[thick] (\xb,\yb)--(\xb+1.5,\yb+1.5)
(\xb+2,\yb)--(\xb+1,\yb)--(\xb+1.5,\yb+1.5)
(\xb+2,\yb)--(\xb+1.5,\yb+1.5)
(\xb+3,\yb)--(\xb+1.5,\yb+1.5)
(\xb+1,\yb)--(\xb+1,\yb+1.5)
;

\end{tikzpicture}
\end{center}
\caption{A graph, which is a $\mathcal{U}$-graph, but not a $\mathcal{U}^{++}\cup\mathcal{U}^{--}$-graph.}\label{ugraph}
\end{figure}

\begin{thm}[Rautenbach and Szwarcfiter \cite{rs}]
For a twin-free graph $G$, the following statements are equivalent.
\begin{itemize}
	\item $G$ is a $\{K_{1,4},K_{1,4}^*,K_{2,3}^*,K_{2,4}^*\}$-free graph. 
	(See Figure \ref{u+-graph} for an illustration.)
	\item $G$ is an almost proper interval graph.
	\item $G$ is a $\mathcal{U}^{++}\cup\mathcal{U}^{--}$-graph.
\end{itemize}
\end{thm} 

Note that an interval representation can assign the same interval to twins and
hence the restriction to twin-free graphs does not weaken the statement but simplifies the description.

The next step is to allow all different types of unit intervals.
The class of $\mathcal{U}$-graphs is a proper superclass of the $\mathcal{U}^{++}\cup\mathcal{U}^{--}$-graphs,
because the graph illustrated in Figure \ref{ugraph} is a $\mathcal{U}$-graph, but not a $\mathcal{U}^{++}\cup\mathcal{U}^{--}$-graph (it contains a $K_{1,4}^*$).
Dourado et al. already made some progress in characterizing this class.

\begin{thm}[Dourado et al. \cite{dlprs}]\label{proper=unit}
For a graph $G$, the following two statements are equivalent.
\begin{itemize}
	\item $G$ is a mixed proper interval graph.
	\item $G$ is a mixed unit interval graph.
\end{itemize}
\end{thm}
They also characterized diamond-free mixed unit interval graphs.
There is another approach by Le and Rautenbach \cite{lr} to understand the class of $\mathcal{U}$-graphs
by restricting the ends of the unit intervals to integers.
They found a infinite list of forbidden induced subgraphs, 
which characterize these so-called \textit{integral} $\mathcal{U}$-\textit{graphs}.



\section{Results}
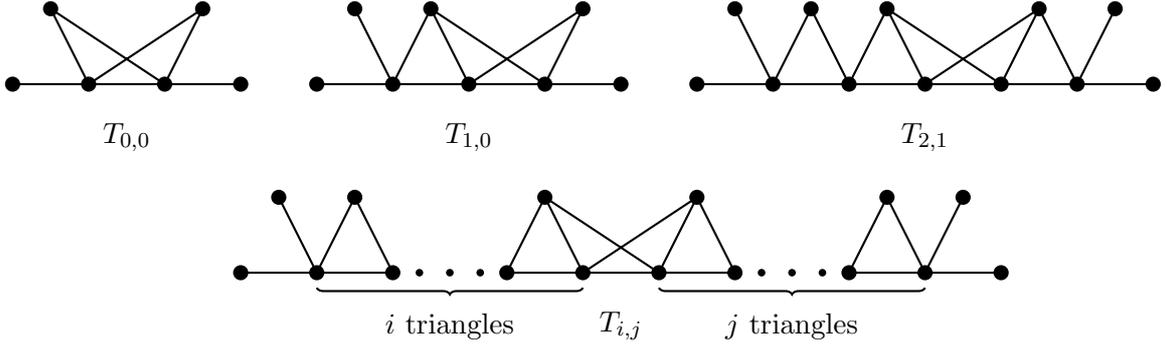
\begin{figure}[t]
\begin{center}
\begin{tikzpicture}[scale=1]
\def\ver{0.1} 
\def\x{1}

\def\xa{0}
\def\ya{0}

\def\xb{4}
\def\yb{0}

\def\xc{9}
\def\yc{0}

\def\xd{3}
\def\yd{-2.5}

\path[fill] (\xa,\ya) circle (\ver);
\path[fill] (\xa+1,\ya) circle (\ver);
\path[fill] (\xa+2,\ya) circle (\ver);
\path[fill] (\xa+3,\ya) circle (\ver);
\path[fill] (\xa+0.5,\ya+1) circle (\ver);
\path[fill] (\xa+2.5,\ya+1) circle (\ver);

\draw[thick] (\xa,\ya)--(\xa+3,\ya)
(\xa+1,\ya)--(\xa+0.5,\ya+1)--(\xa+2,\ya)
(\xa+1,\ya)--(\xa+2.5,\ya+1)--(\xa+2,\ya);

\node (1) at (\xa+1.5,\ya-0.7) {$T_{0,0}$};

\path[fill] (\xb,\yb) circle (\ver);
\path[fill] (\xb+1,\yb) circle (\ver);
\path[fill] (\xb+2,\yb) circle (\ver);
\path[fill] (\xb+3,\yb) circle (\ver);
\path[fill] (\xb+4,\yb) circle (\ver);
\path[fill] (\xb+0.5,\yb+1) circle (\ver);
\path[fill] (\xb+1.5,\yb+1) circle (\ver);
\path[fill] (\xb+3.5,\yb+1) circle (\ver);

\draw[thick] (\xb,\yb)--(\xb+4,\yb)
(\xb+0.5,\yb+1)--(\xb+1,\yb)--(\xb+1.5,\yb+1)--(\xb+2,\yb)--(\xb+3.5,\yb+1)
--(\xb+3,\yb)--(\xb+1.5,\yb+1);

\node (1) at (\xb+2,\yb-0.7) {$T_{1,0}$};

\path[fill] (\xc,\yc) circle (\ver);
\path[fill] (\xc+1,\yc) circle (\ver);
\path[fill] (\xc+2,\yc) circle (\ver);
\path[fill] (\xc+3,\yc) circle (\ver);
\path[fill] (\xc+4,\yc) circle (\ver);
\path[fill] (\xc+5,\yc) circle (\ver);
\path[fill] (\xc+6,\yc) circle (\ver);
\path[fill] (\xc+0.5,\yc+1) circle (\ver);
\path[fill] (\xc+1.5,\yc+1) circle (\ver);
\path[fill] (\xc+2.5,\yc+1) circle (\ver);
\path[fill] (\xc+4.5,\yc+1) circle (\ver);
\path[fill] (\xc+5.5,\yc+1) circle (\ver);

\draw[thick] (\xc,\yc)--(\xc+6,\yc)
(\xc+0.5,\yc+1)--(\xc+1,\yc)--(\xc+1.5,\yc+1)--(\xc+2,\yc)--(\xc+2.5,\yc+1)--(\xc+3,\yc)--(\xc+4.5,\yc+1)
(\xc+5.5,\yc+1)--(\xc+5,\yc)--(\xc+4.5,\yc+1)--(\xc+4,\yc)--(\xc+2.5,\yc+1);

\node (1) at (\xc+3,\yc-0.7) {$T_{2,1}$};

\path[fill] (\xd,\yd) circle (\ver);
\path[fill] (\xd+1,\yd) circle (\ver);
\path[fill] (\xd+2,\yd) circle (\ver);
\path[fill] (\xd+3.5,\yd) circle (\ver);
\path[fill] (\xd+4.5,\yd) circle (\ver);
\path[fill] (\xd+5.5,\yd) circle (\ver);
\path[fill] (\xd+6.5,\yd) circle (\ver);
\path[fill] (\xd+8,\yd) circle (\ver);
\path[fill] (\xd+9,\yd) circle (\ver);
\path[fill] (\xd+10,\yd) circle (\ver);
\path[fill] (\xd+0.5,\yd+1) circle (\ver);
\path[fill] (\xd+1.5,\yd+1) circle (\ver);
\path[fill] (\xd+4,\yd+1) circle (\ver);
\path[fill] (\xd+6,\yd+1) circle (\ver);
\path[fill] (\xd+8.5,\yd+1) circle (\ver);
\path[fill] (\xd+9.5,\yd+1) circle (\ver);

\draw[thick] (\xd,\yd)--(\xd+2,\yd)
(\xd+3.5,\yd)--(\xd+6.5,\yd)
(\xd+8,\yd)--(\xd+10,\yd)
(\xd+0.5,\yd+1)--(\xd+1,\yd)--(\xd+1.5,\yd+1)--(\xd+2,\yd)
(\xd+8,\yd)--(\xd+8.5,\yd+1)--(\xd+9,\yd)--(\xd+9.5,\yd+1)
(\xd+3.5,\yd)--(\xd+4,\yd+1)--(\xd+4.5,\yd)
(\xd+5.5,\yd)--(\xd+6,\yd+1)--(\xd+6.5,\yd)
(\xd+4,\yd+1)--(\xd+5.5,\yd)
(\xd+6,\yd+1)--(\xd+4.5,\yd);

\fill (\xd+2.35,\yd) circle (\ver/2);
\fill (\xd+2.75,\yd) circle (\ver/2);
\fill (\xd+3.15,\yd) circle (\ver/2);

\fill (\xd+6.85,\yd) circle (\ver/2);
\fill (\xd+7.25,\yd) circle (\ver/2);
\fill (\xd+7.65,\yd) circle (\ver/2);

\draw[thick,decoration={brace,mirror,raise=0.2cm},decorate] (\xd+1,\yd) -- (\xd+4.5,\yd)
node [pos=0.5,anchor=north,yshift=-0.4cm] {$i$ triangles};

\draw[thick,decoration={brace,mirror,raise=0.2cm},decorate] (\xd+5.5,\yd) -- (\xd+9,\yd)
node [pos=0.5,anchor=north,yshift=-0.4cm] {$j$ triangles};

\node (1) at (\xd+5,\yd-0.7) {$T_{i,j}$};

\end{tikzpicture}
\end{center}
\caption{The class $\mathcal{T}$.}\label{graphsT}
\end{figure}

In this section we state and prove our main results.
We start by introducing a list of forbidden induced subgraphs.
See Figures \ref{graphsR}, \ref{graphsS}, \ref{graphsS1}, and \ref{graphsT} for illustration.
	Let $\mathcal{R}=\bigcup_{i=0}^{\infty}\{R_i\}$,
$\mathcal{S}=\bigcup_{i=1}^{\infty}\{S_i\}$, $\mathcal{S'}= \bigcup_{i=1}^{\infty}\{S_i'\}$,
and $\mathcal{T}=\bigcup_{i\geq j\geq 0}\{T_{i,j}\}$.
For $k\in\mathbb{N}$ let the graph $Q_k$ arise from the graph $R_k$ by deleting two vertices of degree $1$
that have a common neighbor.
We call the common neighbor of the two deleted vertices and its neighbor of degree $2$ \textit{special vertices} of $Q_k$. 
Note that if a graph $G$ is twin-free, 
then the interval representation of $G$ is injective. 

\begin{figure}[t]
\begin{center}
\begin{tikzpicture}[scale=1]
\def\ver{0.1} 
\def\x{1}

\def\xa{0.5}
\def\ya{0}

\def\xb{4}
\def\yb{0}

\def\xc{8}
\def\yc{0}

\def\xd{3.5}
\def\yd{-2.5}

\path[fill] (\xa+0.5,\ya) circle (\ver);
\path[fill] (\xa+1,\ya+1) circle (\ver);
\path[fill] (\xa+2,\ya+1) circle (\ver);
\path[fill] (\xa+2.5,\ya) circle (\ver);
\path[fill] (\xa+1.5,\ya) circle (\ver);

\draw[thick] (\xa+0.5,\ya)--(\xa+1.5,\ya)--(\xa+1,\ya+1)
(\xa+2,\ya+1)--(\xa+1.5,\ya)--(\xa+2.5,\ya);

\node (1) at (\xa+1.5,\ya-0.7) {$R_0$};

\path[fill] (\xb,\yb) circle (\ver);
\path[fill] (\xb+1,\yb) circle (\ver);
\path[fill] (\xb+2,\yb) circle (\ver);
\path[fill] (\xb+3,\yb) circle (\ver);
\path[fill] (\xb+0.5,\yb+1) circle (\ver);
\path[fill] (\xb+1.5,\yb+1) circle (\ver);
\path[fill] (\xb+2.5,\yb+1) circle (\ver);

\draw[thick] (\xb,\yb)--(\xb+1,\yb)--(\xb+2,\yb)--(\xb+3,\yb)
(\xb+0.5,\yb+1)--(\xb+1,\yb)--(\xb+1.5,\yb+1)--(\xb+2,\yb)--(\xb+2.5,\yb+1);

\node (1) at (\xb+1.5,\yb-0.7) {$R_1$};

\path[fill] (\xc,\yc) circle (\ver);
\path[fill] (\xc+1,\yc) circle (\ver);
\path[fill] (\xc+2,\yc) circle (\ver);
\path[fill] (\xc+3,\yc) circle (\ver);
\path[fill] (\xc+4,\yc) circle (\ver);
\path[fill] (\xc+0.5,\yc+1) circle (\ver);
\path[fill] (\xc+1.5,\yc+1) circle (\ver);
\path[fill] (\xc+2.5,\yc+1) circle (\ver);
\path[fill] (\xc+3.5,\yc+1) circle (\ver);

\draw[thick] (\xc,\yc)--(\xc+1,\yc)--(\xc+2,\yc)--(\xc+3,\yc)--(\xc+4,\yc)
(\xc+0.5,\yc+1)--(\xc+1,\yc)--(\xc+1.5,\yc+1)--(\xc+2,\yc)--(\xc+2.5,\yc+1)--(\xc+3,\yc)--(\xc+3.5,\yc+1);

\node (1) at (\xc+2,\yc-0.7) {$R_2$};

\path[fill] (\xd,\yd) circle (\ver);
\path[fill] (\xd+1,\yd) circle (\ver);
\path[fill] (\xd+2,\yd) circle (\ver);
\path[fill] (\xd+3,\yd) circle (\ver);
\path[fill] (\xd+4.5,\yd) circle (\ver);
\path[fill] (\xd+5.5,\yd) circle (\ver);
\path[fill] (\xd+6.5,\yd) circle (\ver);
\path[fill] (\xd+0.5,\yd+1) circle (\ver);
\path[fill] (\xd+1.5,\yd+1) circle (\ver);
\path[fill] (\xd+2.5,\yd+1) circle (\ver);
\path[fill] (\xd+5,\yd+1) circle (\ver);
\path[fill] (\xd+6,\yd+1) circle (\ver);

\fill (\xd+3.35,\yd) circle (\ver/2);
\fill (\xd+3.75,\yd) circle (\ver/2);
\fill (\xd+4.15,\yd) circle (\ver/2);

\draw[thick] (\xd,\yd)--(\xd+3,\yd)
(\xd+4.5,\yd)--(\xd+6.5,\yd)
(\xd+0.5,\yd+1)--(\xd+1,\yd)--(\xd+1.5,\yd+1)--(\xd+2,\yd)--(\xd+2.5,\yd+1)--(\xd+3,\yd)
(\xd+4.5,\yd)--(\xd+5,\yd+1)--(\xd+5.5,\yd)--(\xd+6,\yd+1);

\draw[thick,decoration={brace,mirror,raise=0.2cm},decorate] (\xd+1,\yd) -- (\xd+5.5,\yd)
node [pos=0.5,anchor=north,yshift=-0.4cm] {$i$ triangles};

\node (1) at (\xd+6,\yd-0.7) {$R_i$};





\end{tikzpicture}
\end{center}
\caption{The class $\mathcal{R}$.}\label{graphsR}
\end{figure}
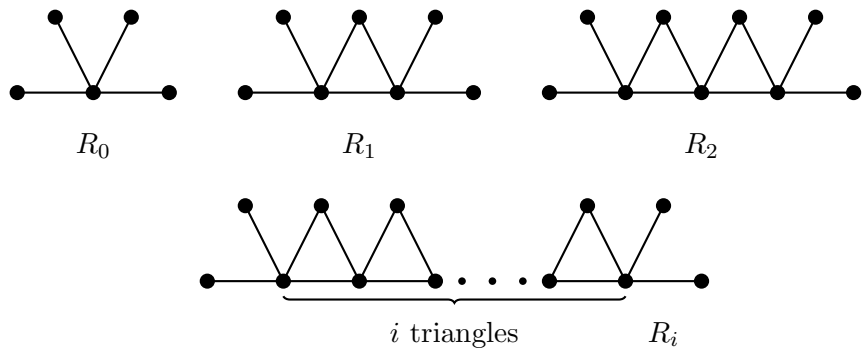


\begin{figure}[t]
\begin{center}
\begin{tikzpicture}[scale=1]
\def\ver{0.1} 
\def\x{1}

\def\xa{0}
\def\ya{0}

\def\xb{0}
\def\yb{0}

\def\xc{4}
\def\yc{0}

\def\xd{9}
\def\yd{0}

\path[fill] (\xb+1,\yb) circle (\ver);
\path[fill] (\xb+2,\yb) circle (\ver);
\path[fill] (\xb+3,\yb) circle (\ver);
\path[fill] (\xb+2,\yb+1.5) circle (\ver);
\path[fill] (\xb+1.5,\yb+1) circle (\ver);
\path[fill] (\xb+2.5,\yb+1) circle (\ver);

\draw[thick] (\xb+1,\yb)--(\xb+2,\yb)--(\xb+3,\yb)
(\xb+1,\yb)--(\xb+1.5,\yb+1)--(\xb+2,\yb)--(\xb+2.5,\yb+1)
--(\xb+2,\yb+1.5)--(\xb+2,\yb)
(\xb+1.5,\yb+1)--(\xb+2,\yb+1.5);

\node (1) at (\xb+2,\yb-0.7) {$S_1$};

\path[fill] (\xc,\yc) circle (\ver);
\path[fill] (\xc+1,\yc) circle (\ver);
\path[fill] (\xc+2,\yc) circle (\ver);
\path[fill] (\xc+3,\yc) circle (\ver);
\path[fill] (\xc+0.5,\yc+1) circle (\ver);
\path[fill] (\xc+1.5,\yc+1) circle (\ver);
\path[fill] (\xc+2.5,\yc+1) circle (\ver);
\path[fill] (\xc+1,\yc+1.5) circle (\ver);

\draw[thick] (\xc,\yc)--(\xc+1,\yc)--(\xc+2,\yc)--(\xc+3,\yc)
(\xc+0.5,\yc+1)--(\xc+1,\yc)--(\xc+1.5,\yc+1)--(\xc+2,\yc)--(\xc+2.5,\yc+1)
(\xc,\yc)--(\xc+0.5,\yc+1)--(\xc+1,\yc+1.5)--(\xc+1.5,\yc+1)
(\xc+1,\yc)--(\xc+1,\yc+1.5);

\node (1) at (\xc+1.5,\yc-0.7) {$S_2$};

\path[fill] (\xd-1,\yd) circle (\ver);
\path[fill] (\xd-0.5,\yd+1) circle (\ver);
\path[fill] (\xd,\yd) circle (\ver);
\path[fill] (\xd+1,\yd) circle (\ver);
\path[fill] (\xd+2,\yd) circle (\ver);
\path[fill] (\xd+3.5,\yd) circle (\ver);
\path[fill] (\xd+4.5,\yd) circle (\ver);
\path[fill] (\xd+5.5,\yd) circle (\ver);
\path[fill] (\xd+0.5,\yd+1) circle (\ver);
\path[fill] (\xd+1.5,\yd+1) circle (\ver);
\path[fill] (\xd+4,\yd+1) circle (\ver);
\path[fill] (\xd+5,\yd+1) circle (\ver);
\path[fill] (\xd,\yd+1.5) circle (\ver);

\fill (\xd+2.35,\yd) circle (\ver/2);
\fill (\xd+2.75,\yd) circle (\ver/2);
\fill (\xd+3.15,\yd) circle (\ver/2);

\draw[thick] (\xd-1,\yd)--(\xd+2,\yd)
(\xd+3.5,\yd)--(\xd+5.5,\yd)
(\xd-1,\yd)--(\xd-0.5,\yd+1)--(\xd,\yd)--(\xd+0.5,\yd+1)--(\xd+1,\yd)--(\xd+1.5,\yd+1)--(\xd+2,\yd)
(\xd+3.5,\yd)--(\xd+4,\yd+1)--(\xd+4.5,\yd)--(\xd+5,\yd+1)
(\xd-0.5,\yd+1)--(\xd,\yd+1.5)--(\xd+0.5,\yd+1)
(\xd,\yd)--(\xd,\yd+1.5);

\draw[thick,decoration={brace,mirror,raise=0.2cm},decorate] (\xd-1,\yd) -- (\xd+4.5,\yd)
node [pos=0.5,anchor=north,yshift=-0.4cm] {$i$ triangles};

\node (1) at (\xd+4.5,\yd-0.7) {$S_i$};

\end{tikzpicture}
\end{center}
\caption{The class $\mathcal{S}$.}\label{graphsS}
\end{figure}



\begin{figure}[t]
\begin{center}
\begin{tikzpicture}[scale=1]
\def\ver{0.1} 
\def\x{1}

\def\xa{4}
\def\ya{0}

\def\xb{0}
\def\yb{0}

\def\xc{8}
\def\yc{0}

\path[fill] (\xa,\ya) circle (\ver);
\path[fill] (\xa+1,\ya) circle (\ver);
\path[fill] (\xa+2,\ya) circle (\ver);
\path[fill] (\xa+3,\ya) circle (\ver);
\path[fill] (\xa+0.5,\ya+1) circle (\ver);
\path[fill] (\xa+1.5,\ya+1) circle (\ver);
\path[fill] (\xa+2.5,\ya+1) circle (\ver);
\path[fill] (\xa+0.5,\ya+0.5) circle (\ver);

\draw[thick] (\xa,\ya)--(\xa+1,\ya)--(\xa+2,\ya)--(\xa+3,\ya)
(\xa+0.5,\ya+1)--(\xa+1,\ya)--(\xa+1.5,\ya+1)--(\xa+2,\ya)--(\xa+2.5,\ya+1)
(\xa,\ya)--(\xa+0.5,\ya+1)--(\xa+0.5,\ya+0.5)--(\xa,\ya)
(\xa+0.5,\ya+0.5)--(\xa+1,\ya);

\node (1) at (\xa+1.5,\ya-0.7) {$S_2'$};

\path[fill] (\xb+1,\yb) circle (\ver);
\path[fill] (\xb+2,\yb) circle (\ver);
\path[fill] (\xb+3,\yb) circle (\ver);
\path[fill] (\xb+1.5,\yb+1) circle (\ver);
\path[fill] (\xb+2.5,\yb+1) circle (\ver);
\path[fill] (\xb+1.5,\yb+0.5) circle (\ver);

\draw[thick] (\xb+1,\yb)--(\xb+2,\yb)--(\xb+3,\yb)
(\xb+1,\yb)--(\xb+1.5,\yb+1)--(\xb+2,\yb)--(\xb+2.5,\yb+1)
(\xb+1.5,\yb+1)--(\xb+1.5,\yb+0.5)--(\xb+1,\yb)
(\xb+1.5,\yb+0.5)--(\xb+2,\yb);

\node (1) at (\xb+2,\yb-0.7) {$S_1'$};

\path[fill] (\xc+1,\yc) circle (\ver);
\path[fill] (\xc+2,\yc) circle (\ver);
\path[fill] (\xc+3,\yc) circle (\ver);
\path[fill] (\xc+4.5,\yc) circle (\ver);
\path[fill] (\xc+5.5,\yc) circle (\ver);
\path[fill] (\xc+6.5,\yc) circle (\ver);
\path[fill] (\xc+1.5,\yc+1) circle (\ver);
\path[fill] (\xc+2.5,\yc+1) circle (\ver);
\path[fill] (\xc+5,\yc+1) circle (\ver);
\path[fill] (\xc+6,\yc+1) circle (\ver);
\path[fill] (\xc+1.5,\yc+0.5) circle (\ver);

\fill (\xc+3.35,\yc) circle (\ver/2);
\fill (\xc+3.75,\yc) circle (\ver/2);
\fill (\xc+4.15,\yc) circle (\ver/2);

\draw[thick] (\xc+1,\yc)--(\xc+3,\yc)
(\xc+4.5,\yc)--(\xc+6.5,\yc)
(\xc+1,\yc)--(\xc+1.5,\yc+1)--(\xc+2,\yc)--(\xc+2.5,\yc+1)--(\xc+3,\yc)
(\xc+4.5,\yc)--(\xc+5,\yc+1)--(\xc+5.5,\yc)--(\xc+6,\yc+1)
(\xc+1.5,\yc+1)--(\xc+1.5,\yc+0.5)--(\xc+1,\yc)
(\xc+1.5,\yc+0.5)--(\xc+2,\yc);

\draw[thick,decoration={brace,mirror,raise=0.2cm},decorate] (\xc+1,\yc) -- (\xc+5.5,\yc)
node [pos=0.5,anchor=north,yshift=-0.4cm] {$i$ triangles};

\node (1) at (\xc+6,\yc-0.7) {$S_i'$};

\end{tikzpicture}
\end{center}
\caption{The class  $\mathcal{S}'$.}\label{graphsS1}
\end{figure}



\begin{lemma}[Dourado et al.\cite{dlprs}]\label{lemma dlprs}
Let $k\in \mathbb{N}$.
\begin{enumerate}[(a)]
	\item Every $\mathcal{U}$-representation of the claw $K_{1,3}$ arises by translation (replacing $I$ by $I+x$ for some $x\in \mathbb{R}$; that is, shifting all intervals by $x$) of the following
  $\mathcal{U}$-representation
	$I: V(K_{1,3})\rightarrow \mathcal{U}$ of $K_{1,3}$,
	where $I(V(K_{1,3}))$ consists of the following intervals
	\begin{itemize}
	\item either $[0,1]$ or $(0,1]$,
	\item $[1,2]$ and $(1,2)$, and
	\item either $[2,3]$ or $[2,3)$.
  \end{itemize}
	\item Every injective $\mathcal{U}$-representation of $Q_k$ arises by translation and inversion (replacing $I$ by $-I$; that is, multiplying all endpoints of the intervals by $-1$)
	of one of the two injective $\mathcal{U}$-representations
	$I: V(Q_k)\rightarrow \mathcal{U}$ of $Q_k$,
	where $I(V(Q_k))$ consists of the following intervals
	\begin{itemize}
	\item either $[0,1]$ or $(0,1]$,
	\item $[1,2]$ and $(1,2)$, and
	\item $[i,i+1]$ and $[i,i+1)$ for $2\leq i \leq k+1$.
  \end{itemize}
	\item The graphs in $\{T_{0,0}\}\cup\mathcal{R}$ are minimal forbidden subgraphs for the class of $\mathcal{U}$-graphs
	with respect to induced subgraphs.
	\item If $G$ is a $\mathcal{U}$-graph, then
	every induced subgraph $H$ in $G$ that is isomorphic to $Q_k$ and every vertex $u^*\in V(G)\setminus V(H)$
	such that $u^*$ is adjacent to exactly one of the two special vertices $x$ of $H$,
	the vertex $u^*$ has exactly one neighbor in $V(H)$, namely $x$. 
\end{enumerate}
\end{lemma}

\begin{lemma}\label{lemma1}
If a graph $G$ is a twin-free mixed unit interval graph,
then $G$ is $\{K_{2,3}^*\}\cup \mathcal{R}\cup \mathcal{S}\cup \mathcal{S'}\cup \mathcal{T}$-free.
\end{lemma}

\noindent
\textit{Proof of Lemma \ref{lemma1}}:
It is easy to see that $G$ is $\{K_{2,3}^*\}$-free.
Lemma \ref{lemma dlprs} (c) shows that $G$ is $\mathcal{R}$-free and
Lemma \ref{lemma dlprs} (d) shows that $G$ is $\mathcal{S}$-free.

Let $k\in \mathbb{N}$. 
Note that the graph $S_k'$ arises from the graph $Q_k$ by adding a vertex $z$ and
joining it to the two special vertices of $Q_k$ and the unique common neighbor of these two vertices.
For contradiction, we assume that $S_k'$ has a $\mathcal{U}$-representation $I$.
By Lemma \ref{lemma dlprs} (b) there are only two possibilities for the $\mathcal{U}$-representation of $Q_k$.
Thus we assume that the subgraph $Q_k$ of $S_k'$ has the interval representation as described in Lemma \ref{lemma dlprs} (b).
In both cases we conclude $\ell(z)=k+1$ and $k+1\in I(z)$.
Thus $r(z)=k+2$ and hence $I(z)\in \{[k+1,k+2],[k+1,k+2)\}$.
Therefore, $G$ is not twin-free, which is a contradiction.
This implies that $G$ is $\mathcal{S}'$-free.

By Lemma \ref{lemma dlprs} (c), $G$ is $T_{0,0}$-free.
Let $C$ be a claw with vertex set $\{c,a_1,a_2,a_3\}$, where $c$ is the center vertex.
Denote by $v_k$ and $w_k$ the special vertices of $Q_k$. 
Note that $T_{k,0}$ arises by the disjoint union of the graph $Q_k$ and $C$,
identifying $v_k$ and $a_1$, and adding the edges $w_kc$ and $v_ka_2$.
For contradiction, we assume that $T_{k,0}$ has a $\mathcal{U}$-representation $I$.
By Lemma \ref{lemma dlprs} (b), we assume without loss of generality that the induced subgraph $Q_k$ of $T_{k,0}$
is represented by exactly the intervals described in Lemma \ref{lemma dlprs} (b).
Thus $I(v_k)=[k+1,k+2]$ and $I(w_k)=[k+1,k+2)$, because $v_ka_2\in E(T_{k,0})$ but $w_ka_2\notin E(T_{k,0})$.
Since $I(v_k)$ is not an open interval and by Lemma \ref{lemma dlprs} (a),
we obtain $I(c)=[k+2,k+3]$
and hence $I(w_k)\cap I(c)=\emptyset$.
This is a contradiction, which implies that $G$ is $\bigcup_{i\geq 0}\{T_{i,0}\}$-free.

Let $i,j\in\mathbb{N}$.
Note that the graph $T_{i,j}$ arises by the disjoint union of $Q_i$ and $Q_j$ and adding three
edges between the special vertices of $Q_i$ and $Q_j$.
We may assume that the intervals of the subgraph $Q_i$ are 
exactly the intervals as described in Lemma \ref{lemma dlprs} (b).
Let $w_i$ (respectively $v_i$) be the vertex of $Q_i$ that has one (two) neighbor(s) in the subgraph $Q_j$;
that is, $I(v_i)=[i+1,i+2]$ and $I(w_i)=[i+1,i+2)$ because $N_{T_{i,j}}(w_i)\subset N_{T_{i,j}}(v_i)$.
Let $w_j$ (respectively $v_j$) be the vertex of $Q_j$ that has one (two) neighbor(s) in the subgraph $Q_i$.
Since the subgraph $Q_j$ has also an interval representation as described in Lemma \ref{lemma dlprs} (b)
and the vertices of $Q_i\setminus\{v_i,w_i\}$ and not joined by an edge to the vertices of $Q_j\setminus\{v_j,w_j\}$,
we conclude that the intervals of the vertices of $Q_j$ arise by an inversion and a translation of the interval representation as described in Lemma \ref{lemma dlprs} (b).
This implies that $I(v_j)=[x,x+1]$ and $I(w_j)=(x,x+1]$ for some $x\in \mathbb{R}$.
Obviously, $x\in[i+1,i+2]$.
If $x=i+2$, then neither $v_i$ is adjacent to $w_j$ nor $w_i$ is adjacent to $v_i$.
If $x\in [i+1,i+2)$, then the intervals of $w_i$ and $w_j$ intersect, which is not possible.
Therefore, $G$ is $\mathcal{T}$-free and this completes the proof.
$\Box$

\bigskip

\noindent
We proceed to our main result.

\begin{thm}\label{mainthm}
A twin-free graph $G$ is a mixed unit interval graph if and only if 
$G$ is a $\{K_{2,3}^*\}\cup\mathcal{R}\cup \mathcal{S}\cup \mathcal{S'}\cup \mathcal{T}$-free interval graph.
\end{thm}

\noindent
\textit{Proof of Theorem \ref{mainthm}}:
We use a similar approach as in \cite{rs}.
By Lemma \ref{lemma1}, we know if $G$ is a twin-free mixed unit interval graph, 
then $G$ is a $\{K_{2,3}^*\}\cup\mathcal{R}\cup \mathcal{S}\cup \mathcal{S'}\cup\mathcal{T}$-free interval graph.
Let $G$ be a twin-free $\{K_{2,3}^*\}\cup\mathcal{R}\cup \mathcal{S}\cup \mathcal{S'}\cup\mathcal{T}$-free interval graph.
We show that $G$ is a mixed proper interval graph.
By Theorem \ref{proper=unit}, this proves Theorem \ref{mainthm}.
Since $G$ is an interval graph,
$G$ has an $\mathcal{I}^{++}$-representation $I$.
As in \cite{rs} we call a pair $(u,v)$ of distinct vertices a \textit{bad pair} if $I(u)\subseteq I(v)$.
Let $I$ be such that the number of bad pairs is as small as possible.
If $I$ has no bad pair, then we are done by Theorem \ref{thmroberts}.
Hence we assume that there is at least one bad pair.
The strategy of the proof is as follows.
Claim \ref{c1} to Claim \ref{c6a} collect properties of $G$ and $I$, 
before we modify our interval representation of $G$ to show that $G$ is a mixed proper interval graph.
In Claim \ref{c7} to Claim 10 we prove that our modification of the interval representation preserves all intersections and non-intersections. 
Claim 1 to Claim \ref{c3} are similar to Claim 1 to Claim 3 in \cite{rs}, respectively.
For the sake of completeness we state the proofs here.

\begin{claim}\label{c1}
If $(u,v)$ is a bad pair,
then there are vertices $x$ and $y$ such that $\ell(v)\leq r(x)< \ell(u)$ and $r(u)<\ell(y)\leq r(v)$.
\end{claim}

\noindent
\textit{Proof of Claim 1:}
For contradiction, we assume the existence of a bad pair $(w,v)$ such that
there is no vertex $x$ with $\ell(v)\leq r(x)<\ell(w)$.
A symmetric argument implies the existence of $y$.
Let $u$ be a vertex such that $\ell(u)$ is as small as possible with respect to $I(u)\subseteq I(v)$.
By our assumption there is no vertex $x$ such that $\ell(v)\leq r(x)<\ell(u)$.
Let $\epsilon$ be the smallest distance between two distinct endpoints of intervals of $I$.
Let $I': V(G)\rightarrow \mathcal{I}^{++}$ be such that
$I'(u)=[\ell(v)- \epsilon/2,r(u)]$,
$I'(v)=[\ell(v),r(v)+\epsilon/2]$,
 and
$I'(z)=I(z)$ for $z\in V(G)\setminus \{u,v\}$.
By the choice of $u$ and $\epsilon$,
we conclude that $I'$ is an interval representation of $G$,
but $I'$ has less bad pairs than $I$,
which is a contradiction to our choice of $I$.
This completes the proof.
$\Box$ 

\bigskip

\noindent
Let $a_1$ and $a_2$ be two distinct vertices.
Claim \ref{c1} implies that $\ell(a_1)\not=\ell(a_2)$
and $r(a_1)\not=r(a_2)$.
Suppose $\ell(a_1)<\ell(a_2)$.
If $r(a_1)=\ell(a_2)$,
then let $\epsilon$ be as in the proof of Claim \ref{c1} and $I': V(G)\rightarrow \mathcal{I}^{++}$ be such that
$I'(a_1)=[\ell(a_1),r(a_1)+\epsilon/2]$,
and
$I'(z)=I(z)$ for $z\in V(G)\setminus \{a_1\}$.
By the choice of $\epsilon$,
we conclude that $I'$ is an interval representation of $G$ with as many bad pairs as $I$.
Therefore, we assume without loss of generality that we chose $I$ such that all endpoints of the intervals of $I$ are distinct.
Hence the inequalities in Claim \ref{c1} are strict inequalities.

\begin{claim}\label{c2}
If $(u,w)$ and $(v,w)$ are bad pairs,
then $u=v$, that is, no interval contains two distinct intervals.
\end{claim}

\noindent
\textit{Proof of Claim \ref{c2}:}
For contradiction,
we assume that there are distinct vertices $u'$, $v'$ and $w$ such that
$(u',w)$ and $(v',w)$ are bad pairs.
Let $u$ be a vertex such that $(u,w)$ is a bad pair
and $\ell(u)$ is as small as possible.
Let $v$ be a vertex such that $(v,w)$ is a bad pair 
and $r(v)$ is as large as possible.
Claim \ref{c1} ensures two distinct vertices $x$ and $y$ such that $\ell(w)< r(x)< \ell(u)$ and $r(v)<\ell(y)< r(w)$.

If $u\not=v$ and $I(u)\cap I(v)=\emptyset$, 
then $G[\{w,x,u,v,y\}]$ is isomorphic to $R_0$, which is a contradiction.
If $u\not=v$ and $I(u)\cap I(v)\not=\emptyset$,
then in the graph $G[\{w,x,u,v,y\}]$ the vertices $u$ and $v$ are twins.
Since $G$ is twin-free, 
$u$ and $v$ do not have the same closed neighborhood in $G$ and hence
there is a vertex $z$, which is adjacent to say $u$ (by symmetry) and not to $v$.
Since $I(u)\subset I(w)$, $z$ is adjacent to $w$.
If $z$ is not adjacent to $x$, then $G[\{w,x,z,v,y\}]$ is isomorphic to $R_0$ and
if $z$ is adjacent to $x$, then $G[\{w,x,z,u,v,y\}]$ is isomorphic to $S_1$, which is a contradiction.

If $u=v$, then there is a vertex $z$ such that $(z,u)$ is a bad pair
because $u'$ or $v'$ is a suitable choice.
We choose $z$ such that $\ell(z)$ is minimal.
Claim \ref{c1} ensures the existence of a vertex $x'$ such that $\ell(u)< r(x')<\ell(z)$.
Note that the choice of $u$ and $z$ guarantees $\ell(x')<\ell(w)$, so $xx'\in E(G)$.
Therefore, $G[\{w,x,x',u,z,y\}]$ is isomorphic to $S_1$, which is a contradiction.
This completes the proof of Claim \ref{c2}. $\Box$

\begin{claim} \label{c3}
If $(u,v)$ and $(u,w)$ are bad pairs,
then $v=w$, that is, no interval is contained in two distinct intervals.
\end{claim}

\noindent
\textit{Proof of Claim \ref{c3}:}
Claim \ref{c2} implies that neither $(v,w)$ nor $(w,v)$ is a bad pair.
Thus we may assume $\ell(w)<\ell(v)<\ell(u)$ and $r(u)<r(w)<r(v)$.
By Claim \ref{c1}, there are vertices $x$ and $y$ such that $\ell(v)< r(x)<\ell(u)$ and $r(u)<\ell(y)< r(w)$.
Now, $G[\{v,w,x,u,y\}]$ is isomorphic to $K_{2,3}^*$, 
which is a contradiction and
completes the proof of Claim \ref{c3}.
$\Box$

\bigskip

\begin{figure}[t]
\begin{center}
\begin{tikzpicture}[scale=1]
\def\ver{0.1} 
\def\x{1}

\def\xa{0}
\def\ya{0}

\def\xb{4}
\def\yb{0}

\def\xc{-4}
\def\yc{-4}

\def\xd{-4}
\def\yd{0}

\path[fill] (\xd,\yd) circle (\ver);
\path[fill] (\xd+1.5,\yd) circle (\ver);
\path[fill] (\xd+3,\yd) circle (\ver);
\path[fill] (\xd+5,\yd) circle (\ver);
\path[fill] (\xd+6.5,\yd) circle (\ver);
\path[fill] (\xd+8,\yd) circle (\ver);
\path[fill] (\xd+10,\yd) circle (\ver);
\path[fill] (\xd+11.5,\yd) circle (\ver);
\path[fill] (\xd+13,\yd) circle (\ver);
\path[fill] (\xd+2.25,\yd+1.5) circle (\ver);
\path[fill] (\xd+5.75,\yd+1.5) circle (\ver);
\path[fill] (\xd+6.5,\yd+1.5) circle (\ver);
\path[fill] (\xd+7.25,\yd+1.5) circle (\ver);
\path[fill] (\xd+10.75,\yd+1.5) circle (\ver);

\draw[thick] (\xd,\yd)--(\xd+3,\yd)
(\xd+5,\yd)--(\xd+8,\yd)
(\xd+10,\yd)--(\xd+13,\yd)
(\xd+1.5,\yd)--(\xd+2.25,\yd+1.5)--(\xd+3,\yd)
(\xd+10,\yd)--(\xd+10.75,\yd+1.5)--(\xd+11.5,\yd)
(\xd+5,\yd)--(\xd+5.75,\yd+1.5)--(\xd+6.5,\yd)--(\xd+7.25,\yd+1.5)-- (\xd+8,\yd)
(\xd+6.5,\yd+1.5)--(\xd+6.5,\yd);

\fill (\xd+3.6,\yd) circle (\ver/2);
\fill (\xd+4,\yd) circle (\ver/2);
\fill (\xd+4.4,\yd) circle (\ver/2);

\fill (\xd+8.6,\yd) circle (\ver/2);
\fill (\xd+9,\yd) circle (\ver/2);
\fill (\xd+9.4,\yd) circle (\ver/2);

\draw[thick,decoration={brace,mirror,raise=1.7cm},decorate] (\xd+1.5,\yd) -- (\xd+6.45,\yd)
node [pos=0.5,anchor=north,yshift=-1.9cm] {$\ell_{u,v}^{\rm max}-1$ triangles};

\draw[thick,decoration={brace,mirror,raise=1.7cm},decorate] (\xd+6.55,\yd) -- (\xd+11.5,\yd)
node [pos=0.5,anchor=north,yshift=-1.9cm] {$r_{u,v}^{\rm max}-1$ triangles};

\node (1) at (\xd+6.5,\yd-0.5) {$v$};
\node (1) at (\xd+6.5,\yd+2) {$u$};
\node (1) at (\xd+5,\yd-0.5) {$x_{u,v}^1$};
\node (1) at (\xd+5.75,\yd+2) {${x_{u,v}^1}'$};

\node (1) at (\xd+1.5,\yd-0.5) {$x_{u,v}^{\ell_{u,v}^{\rm max}-1}$};
\node (1) at (\xd+2.25,\yd+2) {${x_{u,v}^{\ell_{u,v}^{\rm max}-1}}'$};
\node (1) at (\xd,\yd-0.5) {$x_{u,v}^{\ell_{u,v}^{\rm max}}$};

\node (1) at (\xd+8,\yd-0.5) {$y_{u,v}^1$};
\node (1) at (\xd+7.25,\yd+2) {${y_{u,v}^1}'$};

\node (1) at (\xd+11.5,\yd-0.5) {$y_{u,v}^{r_{u,v}^{\rm max}-1}$};
\node (1) at (\xd+10.75,\yd+2) {${y_{u,v}^{r_{u,v}^{\rm max}-1}}'$};
\node (1) at (\xd+13,\yd-0.5) {$y_{u,v}^{r_{u,v}^{\rm max}}$};

\node (1) at (\xd+6.5,\yd-1.3) {$X_{u,v}^0$};
\node (1) at (\xd+5.75,\yd+3.2) {$X_{u,v}^1$};
\node (1) at (\xd+2.5,\yd+3.2) {$X_{u,v}^{\ell_{u,v}^{\rm max}-1}$};
\node (1) at (\xd,\yd+0.8) {$X_{u,v}^{\ell_{u,v}^{\rm max}}$};
\node (1) at (\xd+7.25,\yd+3.2) {$Y_{u,v}^1$};
\node (1) at (\xd+11,\yd+3.2) {$Y_{u,v}^{r_{u,v}^{\rm max}-1}$};
\node (1) at (\xd+13,\yd+0.8) {$Y_{u,v}^{r_{u,v}^{\rm max}}$};

\draw[rotate around={-20:(\xd+5.2,\yd+0.9)}] (\xd+5.2,\yd+0.9) ellipse (1 and 2);
\draw[rotate around={-20:(\xd+1.85,\yd+0.9)}] (\xd+1.85,\yd+0.9) ellipse (1 and 2);
\draw[rotate around={20:(\xd+7.8,\yd+0.9)}] (\xd+7.8,\yd+0.9) ellipse (1 and 2);
\draw[rotate around={20:(\xd+11.0,\yd+0.9)}] (\xd+11.0,\yd+0.9) ellipse (1 and 2);

\draw (\xd-0.2,\yd-0.3) circle (0.7);
\draw (\xd+6.5,\yd-0.3) circle (0.6);
\draw (\xd+13,\yd-0.3) circle (0.7);

\draw[thick] 
(\xc+5.5,\yc)--(\xc+7.5,\yc)
(\xc+6,\yc+0.4)--(\xc+7,\yc+0.4)

(\xc+5.5,\yc+0.2)--(\xc+5.5,\yc-0.2)
(\xc+7.5,\yc+0.2)--(\xc+7.5,\yc-0.2)
(\xc+6,\yc+0.2)--(\xc+6,\yc+0.6)
(\xc+7,\yc+0.2)--(\xc+7,\yc+0.6)

(\xc+4.8,\yc-0.4)--(\xc+5.7,\yc-0.4)
(\xc+4.3,\yc-0.8)--(\xc+5.6,\yc-0.8)

(\xc+4.8,\yc-0.6)--(\xc+4.8,\yc-0.2)
(\xc+5.7,\yc-0.6)--(\xc+5.7,\yc-0.2)
(\xc+4.3,\yc-1)--(\xc+4.3,\yc-0.6)
(\xc+5.6,\yc-1)--(\xc+5.6,\yc-0.6)

(\xc+3.7,\yc+0.4)--(\xc+4.6,\yc+0.4)
(\xc+3.2,\yc)--(\xc+4.5,\yc)

(\xc+3.7,\yc+0.6)--(\xc+3.7,\yc+0.2)
(\xc+4.6,\yc+0.6)--(\xc+4.6,\yc+0.2)
(\xc+3.2,\yc+0.2)--(\xc+3.2,\yc-0.2)
(\xc+4.5,\yc+0.2)--(\xc+4.5,\yc-0.2)

(\xc+2.5,\yc-0.4)--(\xc+3.4,\yc-0.4)
(\xc+2,\yc-0.8)--(\xc+3.3,\yc-0.8)

(\xc+2.5,\yc-0.6)--(\xc+2.5,\yc-0.2)
(\xc+3.4,\yc-0.6)--(\xc+3.4,\yc-0.2)
(\xc+2,\yc-1)--(\xc+2,\yc-0.6)
(\xc+3.3,\yc-1)--(\xc+3.3,\yc-0.6)

(\xc+7.3,\yc-0.4)--(\xc+8.2,\yc-0.4)
(\xc+7.4,\yc-0.8)--(\xc+8.7,\yc-0.8)

(\xc+7.3,\yc-0.6)--(\xc+7.3,\yc-0.2)
(\xc+8.2,\yc-0.6)--(\xc+8.2,\yc-0.2)
(\xc+7.4,\yc-1)--(\xc+7.4,\yc-0.6)
(\xc+8.7,\yc-1)--(\xc+8.7,\yc-0.6)
;

\node (1) at (\xc+6.5,\yc-0.3) {$v$};
\node (1) at (\xc+6.5,\yc+0.7) {$u$};

\node (1) at (\xc+5.1,\yc-1.2) {$X_{u,v}^1$};
\node (1) at (\xc+4.1,\yc+0.9) {$X_{u,v}^2$};
\node (1) at (\xc+2.8,\yc-1.2) {$X_{u,v}^3$};
\node (1) at (\xc+8,\yc-1.2) {$Y_{u,v}^1$};

\fill (\xc+1.3,\yc-0.1) circle (\ver/3);
\fill (\xc+0.9,\yc-0.1) circle (\ver/3);
\fill (\xc+0.5,\yc-0.1) circle (\ver/3);

\fill (\xc+9.4,\yc-0.1) circle (\ver/3);
\fill (\xc+9.8,\yc-0.1) circle (\ver/3);
\fill (\xc+10.2,\yc-0.1) circle (\ver/3);

\end{tikzpicture}
\end{center}
\caption{The structure in $G$ forced by a bad pair $(u,v)$.}\label{structurebadpair}
\end{figure}
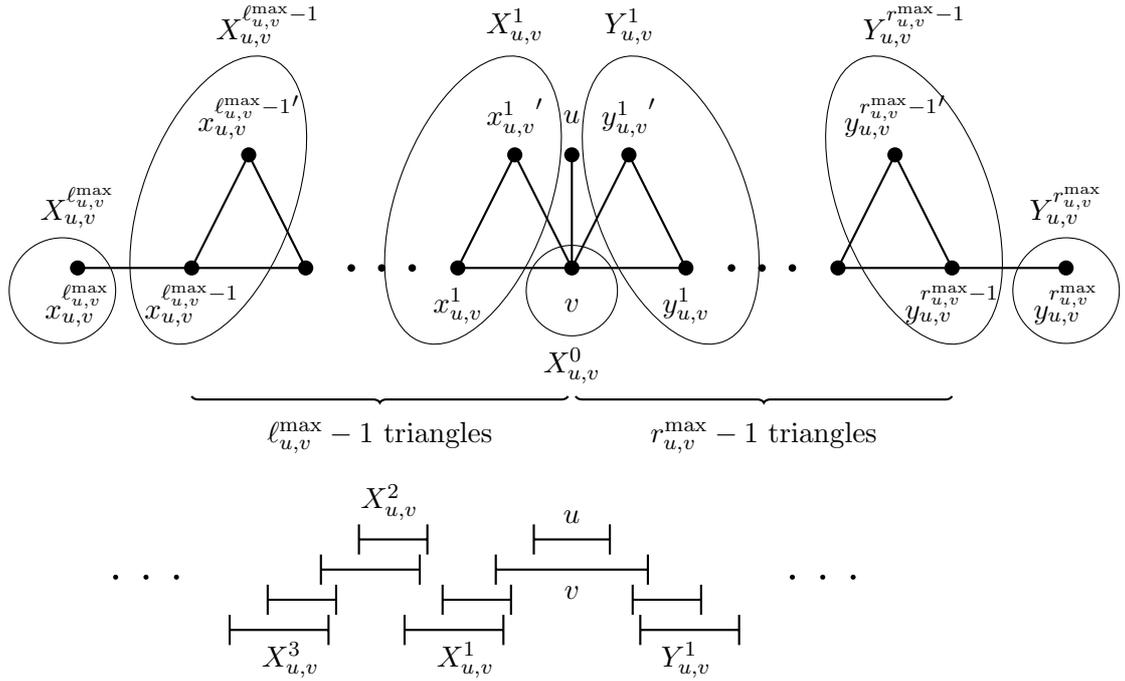

A vertex $x$ is to the \textit{left} (respectively \textit{right}) of a vertex $y$ (in $I$),
if $r(x)<\ell(y)$ (respectively $r(y)<\ell(x)$).
Two adjacent vertices $x$ and $y$ are \textit{distinguishable} by vertices to the left (respectively right) of them,
if there is a vertex $z$, which is adjacent to exactly one of them and to the left (respectively right) of one of them.
The vertex $z$ \textit{distinguishes} $x$ and $y$.
Next, we show that for a bad pair $(u,v)$ there is the structure as shown in Figure \ref{structurebadpair} in $G$.
We introduce a positive integer $\ell_{u,v}^{\rm max}$ that, roughly speaking, indicates how large this structure is.

For a bad pair $(u,v)$ let $v=X_{u,v}^0$ and let $X_{u,v}^1$ be the set of vertices that are adjacent to $v$ and
to the left of $u$. 
Let $y_{u,v}$ be a vertex to the right of $u$ and adjacent to $v$.
Claim \ref{c1} guarantees $|X_{u,v}^1|\geq 1$ and the existence of $y_{u,v}$.
If $|X_{u,v}^1|=1$, 
then let $\ell_{u,v}^{\rm max}=1$ and we stop here.
Suppose $|X_{u,v}^1|\geq 2$.
Since $G$ is $R_0$-free, $X_{u,v}^1$ is a clique and
since $G$ is $S_1'$-free, we conclude $|X_{u,v}^1|=2$.
Let $\{x,x'\}= X_{u,v}^1$ such that $r(x)< r(x')$.
For contradiction, we assume that there is a vertex $z$ to the right of $x$ that distinguishes $x$ and $x'$.
We conclude $\ell(v)<\ell(z)$.
By Claim \ref{c2}, $r(v)< r(z)$.
This implies that $(u,z)$ is a bad pair, which contradicts Claim \ref{c3}.
Thus $z$ does not exist.
In addition $(x,x')$ is not a bad pair, 
otherwise Claim \ref{c1} guarantees a vertex $z$ such that $r(x)<\ell(z)<r(x')$, which is a contradiction.
Thus $\ell(x)<\ell(x')<r(x)< r(x')$.
Let $x_{u,v}^1=x$ and ${x_{u,v}^1}'=x'$.
Note that $N_G({x^1_{u,v}}')\subset N_G(x^1_{u,v})$.

Let $X_{u,v}^2=N_G(x^1_{u,v})\setminus N_G({x^1_{u,v}}')$.
Note that all vertices in $X_{u,v}^2$ are to the left of ${x^1_{u,v}}'$.
Since $G$ is twin-free, $|X_{u,v}^2|\geq 1$.
If $|X_{u,v}^2|=1$,
then let $\ell_{u,v}^{\rm max}=2$ and we stop here.
Suppose $|X_{u,v}^2|\geq 2$.
Since $G$ is $R_1$-free, $X_{u,v}^2$ is a clique and
since $G$ is $S_2'$-free, we conclude $|X_{u,v}^2|=2$.
Let $\{x,x'\}= X_{u,v}^2$ such that $r(x)< r(x')$.
For contradiction, we assume that there is a vertex $z$ to the right of $x$ that distinguishes $x$ and $x'$.
Since $z\notin X_{u,v}^2$,
we conclude $\ell({x_{u,v}^1}')< r(z)$.
If $r(z)<\ell(v)$, then $G[\{z,x,x',{x_{u,v}^1},{x_{u,v}^1}',v,u,y_{u,v}\}]$ is isomorphic to $S_2$, which is a contradiction.
Thus $\ell(v)< r(z)$.
If $r(z)< \ell(u)$, then $|X_{u,v}^1|=3$, which is a contradiction.
Thus $\ell(u)< r(z)$.
If $r(u)< r(z)$, then $(u,v)$ and $(u,z)$ are bad pairs, which is a contradiction to Claim \ref{c3}.
Thus $\ell(u)< r(z)<r(u)$.
Now $G[\{z,x',{x_{u,v}^1}',v,u,y_{u,v}\}]$ is isomorphic to $T_{0,0}$, which is the final contradiction.

Note that $(x,x')$ is not a bad pair, 
otherwise Claim 1 guarantees a vertex $z$ such that $r(x)<\ell(z)< r(x')$, which is a contradiction.
Thus $\ell(x)<\ell(x')<r(x)< r(x')$.
Let $x_{u,v}^2=x$ and ${x_{u,v}^2}'=x'$.
Note that $N_G({x^2_{u,v}}')\subset N_G(x^2_{u,v})$.
Let $X_{u,v}^3=N_G(x^2_{u,v})\setminus N_G({x^2_{u,v}}')$.
Note that all vertices in $X_{u,v}^3$ are to the left of ${x^2_{u,v}}'$.

\noindent
We assume that for $k\geq 3$, $i\in[k-1]$ and $j\in[k]$
\begin{itemize}
	\item we defined $X_{u,v}^j$,
	\item $|X_{u,v}^i|=2$ holds,
	\item we defined $x_{u,v}^i$ and ${x_{u,v}^i}'$,
	\item $\ell(x_{u,v}^i)<\ell({x_{u,v}^i}')<r(x_{u,v}^i)<r({x_{u,v}^i}')$ holds,
	\item the vertices in $X_{u,v}^{i+1}$ are to the left of ${x^{i}_{u,v}}'$, and
	\item the vertices in $X_{u,v}^i$ are not distinguishable to the right.
\end{itemize}

If $|X_{u,v}^k|=1$,
then let $\ell_{u,v}^{\rm max}=k$ and we stop here.
Suppose $|X_{u,v}^k|\geq 2$.
Since $G$ is $R_{k-1}$-free, $X_{u,v}^k$ is a clique and
since $G$ is $S_k'$-free, we obtain $|X_{u,v}^k|=2$.
Let $\{x,x'\}= X_{u,v}^k$ such that $r(x)< r(x')$.
For contradiction, we assume that there is a vertex $z$ to the right of $x$ that distinguishes $x$ and $x'$.
Since $z\notin X_{u,v}^k$,
we conclude $\ell({x_{u,v}^{k-1}}')< r(z)$.
If $r(z)<\ell(x_{u,v}^{k-2})$, then $G[\{z,x,x',v,u,y_{u,v}\}\cup \bigcup_{i=1}^{k-1}X_{u,v}^i]$ is isomorphic to $S_k$, which is a contradiction.
Thus $\ell(x_{u,v}^{k-2})< r(z)$.
If $r(z)<\ell({x_{u,v}^{k-2}}')$, then $|X_{u,v}^{k-1}|=3$, which is a contradiction.
Thus $\ell({x_{u,v}^{k-2}}')< r(z)$.
If $r(z)<\ell(x_{u,v}^{k-3})$, then $G[\{z,x',{x_{u,v}^{k-1}}',v,u,y_{u,v}\}\cup \bigcup_{i=1}^{k-2}X_{u,v}^i]$ is isomorphic to $T_{k-3,0}$, which is a contradiction.
Thus $\ell(x_{u,v}^{k-3})< r(z)$.
If $r(z)< r(x_{u,v}^{k-2})$, then $|X_{u,v}^{k-2}|=3$, which is a contradiction.
Thus $r(x_{u,v}^{k-2})<r(z)$ and hence $({x_{u,v}^{k-1}}',z)$ and $(x_{u,v}^{k-2},z)$ are bad pairs, which is a contradiction to Claim \ref{c2}.
Thus $x,x'$ are not distinguishable to the right.
We obtain that $(x,x')$ is not a bad pair, 
otherwise Claim 1 guarantees a vertex $z$ such that $r(x)<\ell(z)< r(x')$, which is a contradiction.
Thus $\ell(x)<\ell(x')<r(x)< r(x')$.
Let $x_{u,v}^k=x$ and ${x_{u,v}^k}'=x'$.
Note that $N_G({x^k_{u,v}}')\subset N_G(x^k_{u,v})$.
Let $X_{u,v}^{k+1}=N_G(x^k_{u,v})\setminus N_G({x^k_{u,v}}')$.
Note that all vertices in $X_{u,v}^{k+1}$ are to the left of ${x^k_{u,v}}'$.

\noindent
By induction this leads to the following properties.

\begin{claim} \label{c4}
If $(u,v)$ is a bad pair, $k\in [\ell_{u,v}^{\rm max}-1]$,
then the following holds:
\begin{enumerate}[(a)]
	\item $|X_{u,v}^k|= 2$.
	\item The vertices in $X_{u,v}^k$ are not distinguishable by vertices to the right of them.
	\item We have $\ell(x_{u,v}^i)<\ell({x_{u,v}^i}')<r(x_{u,v}^i)<r({x_{u,v}^i}')$, 
	that is $(x^k_{u,v},{x^k_{u,v}}')$ and $({x^k_{u,v}}',x^k_{u,v})$ are not bad pairs.
\end{enumerate}
\end{claim}

\noindent
Note that $\ell_{u,v}^{\rm max}$ is the smallest integer $k$ such that $|X_{u,v}^{k-1}|= 2$ and
$|X_{u,v}^k|=1$.

\begin{claim} \label{c5}
If $(u,v)$ is a bad pair and $k\in [\ell_{u,v}^{\rm max}-1]$,
then the following holds.
\begin{enumerate}[(a)]
	\item ${x_{u,v}^k}'$ is not contained in a bad pair.
	\item There is no vertex $z\in V(G)$ such that $(x_{u,v}^k,z)$ is a bad pair.
\end{enumerate}
\end{claim}

\noindent
\textit{Proof of Claim \ref{c5}:}
(a): For contradiction, we assume that there is a vertex $z\in V(G)$ such that $({x_{u,v}^k}',z)$ is a bad pair.
Trivially $z\notin \{\{v,y,u\} \cup \bigcup_{i=1}^{\ell_{u,v}^{\rm max}} X_{u,v}^i \}$.
We have $r({x_{u,v}^k}')<r(z)$ and $\ell(z)<\ell({x_{u,v}^k}')$.
In addition $\ell(x_{u,v}^k)<\ell(z)$, otherwise $(x_{u,v}^k,z)$ is also a bad pair, which contradicts Claim \ref{c2}.
Claim \ref{c1} implies the existence of a vertex $a$, such that $\ell(z)< r(a)< \ell({x_{u,v}^k}')$.

Let $k=1$.
If $r(z)<\ell(u)$, then $z\in X_{u,v}^1$, which is a contradiction to $|X_{u,v}^1|=2$.
Thus $\ell(u)< r(z)$.
If $r(z)< r(u)$, then $G[\{a,z,{x_{u,v}^k}',u,v,y\}]$ is isomorphic to $T_{0,0}$, which is a contradiction.
Thus $r(u)<r(z)$ and now $(u,z)$ is a bad pair, which is a contradiction to Claim \ref{c2}.

Let $k\geq2$.
If $r(z)<\ell({x_{u,v}^{k-1}}')$, then $z\in X_{u,v}^k$, which is a contradiction to $|X_{u,v}^k|=2$.
Thus $\ell({x_{u,v}^{k-1}}')< r(z)$.
If $r(z)< \ell(x_{u,v}^{k-2})$, then $G[\{a,z,{x_{u,v}^{k}}',v,u,y\}\cup \bigcup_{i=1}^{k-1}X_{u,v}^i]$ is isomorphic to $T_{k-1,0}$.
Thus $\ell(x_{u,v}^{k-2})< r(z)$.
If $r(z)< r(x_{u,v}^{k-1})$, then $z\in X_{u,v}^{k-1}$, which is a contradiction to $|X_{u,v}^{k-1}|=2$.
Thus $r(x_{u,v}^{k-1})<r(z)$, but now $(x_{u,v}^{k-1},z)$ is also a bad pair, 
which is a contradiction to Claim \ref{c2} and completes this part of the proof.

For contradiction, we assume that there is a vertex $z\in V(G)$ such that $(z,{x_{u,v}^k}')$ is a bad pair.
By Claim \ref{c1}, $\ell({x_{u,v}^k}')<\ell(z)$ and $r(z)<r({x_{u,v}^k}')$.
By Claim \ref{c3}, $r(x_{u,v}^k)<r(z)$.
Let $y_z$ be the vertex guaranteed by Claim \ref{c1} such that $r(z)<\ell(y_z)$,
but this contradicts Claim \ref{c4} (b).

(b): For contradiction, we assume the existence of a vertex $z\in V(G)$
such that $(x_{u,v}^k,z)$ is a bad pair.
Trivially $z\not= {x_{u,v}^k}'$.
If $r(z)< r({x_{u,v}^k}')$, then this contradicts Claim \ref{c4} (a), that is $|X^k_{u,v}|=2$ and
if $r({x_{u,v}^k}')<r(z)$, then $({x_{u,v}^k}',z)$ is also a bad pair and this contradicts Claim \ref{c2}.
This completes the proof of Claim \ref{c5}.
$\Box$

\bigskip

For a bad pair $(u,v)$ define $Y_{u,v}^k$ as $X_{u,v}^k$ by interchanging in the definition right by left.
Let $r_{u,v}^{\rm max}$ be the smallest integer $k$ such that $|Y_{u,v}^{k-1}|= 2$ and $|Y_{u,v}^{k}|= 1$.
By symmetry, one can prove a ``y''-version of Claim~\ref{c4}, Claim~\ref{c5} and Claim~\ref{c6a} (a) and (b).
Let $\{y_{u,v}^k,{y_{u,v}^k}'\}=Y_{u,v}^k$ such that $N_G({y^k_{u,v}}')\subset N_G(y^k_{u,v})$ for $k\leq r_{u,v}^{\rm max}-1$.

\begin{claim}\label{c6a}
Let $(u,v)$ and $(w,z)$ be bad pairs and $k\in[\ell_{u,v}^{\rm max}]$.
\begin{enumerate}[(a)]
	\item If $X_{u,v}^{k}\cap X_{w,z}^{\tilde{k}}\not=\emptyset$,
	then $x_{u,v}^{k-1}=x_{w,z}^{\tilde{k}-1}$ for $\tilde{k}\in [\ell_{w,z}^{\rm max}]$.
	\item If $X_{u,v}^{k}\cap X_{w,z}^{\tilde{k}}\not=\emptyset$, 
	then $X_{u,v}^{k}=X_{w,z}^{\tilde{k}}$ for $\tilde{k}\in [\ell_{w,z}^{\rm max}]$.
	\item If $X_{u,v}^{k}\cap Y_{w,z}^{\tilde{k}}\not=\emptyset$,
	then $X_{u,v}^{k}\cap Y_{w,z}^{\tilde{k}}=x_{u,v}^{k}= y_{w,z}^{\tilde{k}}$ for $\tilde{k}\in [r_{w,z}^{\rm max}]$
\end{enumerate}
\end{claim}

\noindent
\textit{Proof of Claim \ref{c6a}:}
(a): For contradiction we assume $x_{u,v}^{k-1}\not=x_{w,z}^{\tilde{k}-1}$.
Without loss of generality we assume $\ell(x_{u,v}^{k-1})<\ell(x_{w,z}^{\tilde{k}-1})$.
Note that $x_{w,z}^{\tilde{k}-1}$ is adjacent to the vertices in $X_{u,v}^{k}\cap X_{w,z}^{\tilde{k}}$.
Since the vertices in $X_{u,v}^{k}$ are not distinguishable to the right,
we conclude $\ell(x_{w,z}^{\tilde{k}-1})<r(x_{u,v}^{k})$.

First, we suppose $k=1$.
Thus $v=x_{u,v}^{k-1}$.
If $r(x_{w,z}^{\tilde{k}-1})< r(v)$, then $(x_{w,z}^{\tilde{k}-1},v)$ is a bad pair and this contradicts Claim \ref{c2}
and if $r(x_{w,z}^{\tilde{k}-1})> r(v)$, then $(u,x_{w,z}^{\tilde{k}-1})$ is a bad pair and this contradicts Claim \ref{c3}.
Now we suppose $k\geq 2$.
If $r({x_{u,v}^{k-1}}')< r(x_{w,z}^{\tilde{k}-1})$,
then $({x_{u,v}^{k-1}}',x_{w,z}^{\tilde{k}-1})$ is a bad pair,
which contradicts Claim \ref{c5} (a).
Thus $r(x_{w,z}^{\tilde{k}-1})< r({x_{u,v}^{k-1}}')$.
If $r({x_{u,v}^{k-1}})<r(x_{w,z}^{\tilde{k}-1})$,
then $x_{w,z}^{\tilde{k}-1}\in X_{u,v}^{k-1}$,
which implies $|X_{u,v}^{k-1}|=3$ and hence contradicts Claim \ref{c4} (a).
Thus $r(x_{w,z}^{\tilde{k}-1})<r({x_{u,v}^{k-1}})$.
Therefore, $(x_{w,z}^{\tilde{k}-1},x_{u,v}^{k-1})$ is a bad pair.
Claim \ref{c1} implies the existence of a vertex $a$ which is to the left of $x_{w,z}^{\tilde{k}-1}$
and adjacent to $x_{u,v}^{k-1}$.
Thus $a\in X_{u,v}^{k}$.
However, $r(a)<r(x_{u,v}^{k})$, which contradicts Claim \ref{c4} (c).
This is the final contradiction and this completes the proof of Claim \ref{c6a} (a).

(b): If $|X_{u,v}^{k}|=|X_{w,z}^{\tilde{k}}|=1$,
then there is nothing to show.
Thus we assume, $|X_{u,v}^{k}|=2$.
Note that by Claim \ref{c6a} (a), 
$x_{u,v}^{k-1}=x_{w,z}^{\tilde{k}-1}$.
If ${x_{u,v}^{k}}'\in X_{w,z}^{\tilde{k}}$,
then ${x_{u,v}^{k}}\in X_{w,z}^{\tilde{k}}$ and we are done.
Thus we assume ${x_{u,v}^{k}}'\notin X_{w,z}^{\tilde{k}}$.
Since $X_{u,v}^{k}\cap X_{w,z}^{\tilde{k}}\not=\emptyset$,
we conclude ${x_{u,v}^{k}}\in X_{w,z}^{\tilde{k}}$.
Hence $w$ or ${x_{w,z}^{\tilde{k}-1}}'$ distinguishes the vertices in $X_{u,v}^{k}$ to the right of them,
which is a contradiction to Claim \ref{c4} (b).
This completes the proof.

(c): If $|X_{u,v}^{k}|=|Y_{w,z}^{\tilde{k}}|=1$,
then there is nothing to show.
Thus we assume by symmetry $|Y_{w,z}^{\tilde{k}}|=2$.
First, we assume for contradiction ${y_{w,z}^{\tilde{k}}}'=x_{u,v}^{k}$. 
Note that $\ell(y_{w,z}^{\tilde{k}-1})<\ell({y_{w,z}^{\tilde{k}}}')$ and $r(y_{w,z}^{\tilde{k}-1})<r({y_{w,z}^{\tilde{k}}}')$.

Suppose $|X_{u,v}^{k}|=1$.
If $\ell(x_{u,v}^{k-1})<r(y_{w,z}^{\tilde{k}-1})$,
then $y_{w,z}^{\tilde{k}-1}\in X_{u,v}^k$, which is a contradiction to $|X_{u,v}^k|=1$.
Thus $r(y_{w,z}^{\tilde{k}-1})<\ell(x_{u,v}^{k-1})$.
Note that $\ell({y_{w,z}^{\tilde{k}}}')<\ell({y_{w,z}^{\tilde{k}}})<r({y_{w,z}^{\tilde{k}-1}})$ and
$r({y_{w,z}^{\tilde{k}}}')<r({y_{w,z}^{\tilde{k}}})$.
Suppose $k=1$.
If $r(y_{w,z}^{\tilde{k}})<\ell(u)$, then $y_{w,z}^{\tilde{k}}\in X_{u,v}^k$, which is a contradiction to $|X_{u,v}^k|=1$.
If $\ell(u)<r(y_{w,z}^{\tilde{k}})<r(u)$, then $G[\{x_{w,z}^1,w,z,u,v,y_{u,v}^1\}\cup \bigcup_{i=1}^{\tilde{k}}Y_{w,z}^i]$ is isomorphic to $T_{\tilde{k},0}$, which is a contradiction.
If $r(u)<r(y_{w,z}^{\tilde{k}})$, then $(u,y_{w,z}^{\tilde{k}})$ is a bad pair, which is a contradiction to Claim \ref{c3}.
Now we suppose $k\geq 2$.
If $r(y_{w,z}^{\tilde{k}})<\ell({x_{u,v}^{k-1}}')$, 
then $y_{w,z}^{\tilde{k}}\in X_{u,v}^k$, which is a contradiction to $|X_{u,v}^k|=1$.
If $\ell({x_{u,v}^{k-1}}')<r(y_{w,z}^{\tilde{k}})<\ell(x_{u,v}^{k-2})$, 
then $G[\{x_{w,z}^1,w,z,u,v,y_{u,v}^1\} \cup \bigcup_{i=1}^{\tilde{k}}Y_{w,z}^i\cup \bigcup_{i=1}^{k-1}X_{u,v}^i]$ is isomorphic to $T_{\tilde{k},k-1}$, which is a contradiction.
If $\ell(x_{u,v}^{k-2})<r(y_{w,z}^{\tilde{k}})<\ell({x_{u,v}^{k-1}}')$,
then $y_{w,z}^{\tilde{k}}\in X_{u,v}^{k-1}$ and hence $|X_{u,v}^{k-1}|=3$, which is a contradiction to Claim \ref{c4} (a).
If $\ell({x_{u,v}^{k-1}}')<r(y_{w,z}^{\tilde{k}})$, 
then $({x_{u,v}^{k-1}}',y_{w,z}^{\tilde{k}})$ is a bad pair, which is a contradiction to Claim \ref{c5} (a).

This shows $|X_{u,v}^{k}|\not=1$ and thus we suppose $|X_{u,v}^{k}|=2$.
If $\ell(x_{u,v}^{k-1})<r(y_{w,z}^{\tilde{k}-1})$,
then $y_{w,z}^{\tilde{k}-1}\in X_{u,v}^k$, which is a contradiction to $|X_{u,v}^k|=2$.
Thus $r(y_{w,z}^{\tilde{k}-1})<\ell(x_{u,v}^{k-1})$.
Note that $\ell({y_{w,z}^{\tilde{k}}}')<\ell({y_{w,z}^{\tilde{k}}})<r({y_{w,z}^{\tilde{k}-1}})$ and
$r({y_{w,z}^{\tilde{k}}}')<r({y_{w,z}^{\tilde{k}}})$.
If $\ell({x_{u,v}^k}')<r(y_{w,z}^{\tilde{k}-1})$,
then ${x_{u,v}^k}'=y_{w,z}^{\tilde{k}}$.
Thus $\{{x_{u,v}^k}',x_{u,v}^k\}=Y_{w,z}^{\tilde{k}}$.
By Claim \ref{c4} (b), these vertices are not distinguishable to the right and to the left.
Thus they are twins, which is a contradiction.
Thus $r(y_{w,z}^{\tilde{k}-1})<\ell({x_{u,v}^k}')$.
Note that $\ell({y_{w,z}^{\tilde{k}}}')<\ell({y_{w,z}^{\tilde{k}}})<r({y_{w,z}^{\tilde{k}-1}})$.
If $r(y_{w,z}^{\tilde{k}})<r({x_{u,v}^k}')$, 
then $y_{w,z}^{\tilde{k}}\in  X_{u,v}^k$,
which is a contradiction to $| X_{u,v}^k|=2$ and
if $r({x_{u,v}^k}')<r(y_{w,z}^{\tilde{k}})$,
then $({x_{u,v}^k}',y_{w,z}^{\tilde{k}})$ is a bad pair,
which is a contradiction to Claim \ref{c5} (a).
This shows ${y_{w,z}^{\tilde{k}}}'\not=x_{u,v}^{k}$.
A totally symmetric argumentation shows $y_{w,z}^{\tilde{k}}\not={x_{u,v}^{k}}'$.

To complete the proof, we show that ${y_{w,z}^{\tilde{k}}}'\not={x_{u,v}^{k}}'$.
For contradiction, we assume ${y_{w,z}^{\tilde{k}}}'={x_{u,v}^{k}}'$.
If $\ell(x_{u,v}^{k-1})<r(y_{w,z}^{\tilde{k}-1})$,
then $x_{u,v}^{k-1}=y_{w,z}^{\tilde{k}}$.
Thus $G[\{x_{w,z}^1,w,z,u,v,y_{u,v}^1\}\cup \bigcup_{i=1}^{\tilde{k}}Y_{w,z}^i\cup \bigcup_{i=1}^{k-1}X_{u,v}^i]$ is isomorphic to $R_{k+\tilde{k}-1}$, 
which is a contradiction.
Hence we assume $r(y_{w,z}^{\tilde{k}-1})<\ell(x_{u,v}^{k-1})$.
If $\ell(x_{u,v}^k)<\ell(y_{w,z}^{\tilde{k}-1})$,
then $(y_{w,z}^{\tilde{k}-1},x_{u,v}^k)$ is a bad pair, which is a contradiction to the ``y''-version of Claim \ref{c5} (b).
Hence we assume $\ell(y_{w,z}^{\tilde{k}-1})<\ell(x_{u,v}^k)$.
If $x_{u,v}^k\in Y_{w,z}^{\tilde{k}}$,
then this is a contradiction to the ``y''-version of Claim \ref{c4} (a), because  $\ell(x_{u,v}^k)<\ell({y_{w,z}^{\tilde{k}}}')$.
Suppose $\tilde{k}=1$.
Since $x_{u,v}^k\notin Y_{w,z}^{\tilde{k}}$,
we conclude $x_{u,v}^kw\in E(G)$.
If $\ell(w)<\ell(x_{u,v}^k)$,
then $G[\{x_{w,z}^1,w,z,u,v,y_{u,v}^1\}\cup \bigcup_{i=1}^{k}X_{u,v}^i ]$ is isomorphic to $T_{k,0}$, 
which is a contradiction.
If $\ell(x_{u,v}^k)<\ell(w)$, 
then $(w,x_{u,v}^k)$ is a bad pair, which is a contradiction to Claim \ref{c3}.
Hence we suppose $\tilde{k}\geq2$.
Note that $\ell(x^k_{u,v})<r({y_{w,z}^{\tilde{k}-1}}')$.
If $r(y_{w,z}^{\tilde{k}-2})<\ell(x_{u,v}^k)$,
then $G[\{x_{w,z}^1,w,z,u,v,y_{u,v}^1\}\cup \bigcup_{i=1}^{\tilde{k}-1}Y_{w,z}^i\cup \bigcup_{i=1}^{k}X_{u,v}^i]$ is isomorphic to $T_{\tilde{k}-1,k}$.
If $\ell({y_{w,z}^{\tilde{k}-1}}')<\ell(x_{u,v}^k)<r(y_{w,z}^{\tilde{k}-2})$,
then $x_{u,v}^k \in Y_{w,z}^{\tilde{k}-1}$, which is a contradiction to the ``y''-version of Claim \ref{c4} (a).
If $\ell(x_{u,v}^k)<\ell({y_{w,z}^{\tilde{k}-1}}')$,
then $({y_{w,z}^{\tilde{k}-1}}',x_{u,v}^k)$ is a bad pair, 
which is a contradiction to the ``y''-version of Claim \ref{c5} (a).
This completes the proof of Claim \ref{c6a}.
$\Box$

\bigskip

Next, we define step by step new interval representations of $G$ as follows.
First we shorten the intervals of $X_{u,v}^{k}$ for every bad pair $(u,v)$ and $k\in [\ell_{u,v}^{\rm max}]$.
Let $I':V(G)\rightarrow \mathcal{I}^{++}$ be such that $I'(x)=[\ell(x),\ell(x_{u,v}^{k-1})]$ if $x\in X_{u,v}^k$ for some bad pair $(u,v)$
and $I'(x)=I(x)$ otherwise.
By Claim~\ref{c6a}~(a), $I'$ is well-defined; that is,
if $x\in X_{u,v}^k\cap X_{w,z}^{\tilde{k}}$, then $\ell(x_{u,v}^{k-1})=\ell(x_{w,z}^{\tilde{k}-1})$.
Let $\ell'(x)$ and $r'(x)$ be the left and right endpoint of the interval $I'(x)$ for $x\in V(G)$, respectively.

\begin{claim}\label{c7} 
$I'$ is an interval representation of $G$.
\end{claim}

\noindent
\textit{Proof of Claim \ref{c7}:}
Trivially, if two intervals do not intersect in $I$, then they do not intersect in $I'$.
For contradiction, we assume that there are two vertices $a,b\in V(G)$ such that
$I(a)\cap I(b)\not= \emptyset$ and $I'(a)\cap I'(b)= \emptyset$.
At least one interval is shortened by changing the interval representation.
Say $a\in X_{u,v}^k$ for some bad pair $(u,v)$ and $k\in [\ell_{u,v}^{\rm max}]$.
Hence $b\not=x_{u,v}^{k-1}$ and $\ell(x_{u,v}^{k-1})<\ell(b)$ and by Claim~\ref{c4}~(b),
$\ell(b)< r({x_{u,v}^{k}})$.
We conclude that $(b,x_{u,v}^{k-1})$ is not a bad pair, 
otherwise Claim \ref{c1} implies the existence of a vertex $z\in X_{u,v}^{k}$ to the left of $b$,
but $z\notin\{x_{u,v}^{k},{x_{u,v}^{k}}'\}$, which is a contradiction to Claim \ref{c4} (a).
Thus $r(x_{u,v}^{k-1})<r(b)$.
If $k=1$, then $(u,b)$ is also a bad pair, which is a contradiction to Claim \ref{c3}.
Thus $k\geq 2$.
Since $\ell(b)< r({x_{u,v}^{k}})$,
we obtain $\ell(b)<\ell({x_{u,v}^{k-1}}')$.
Since $({x_{u,v}^{k-1}}',b)$ is not a bad pair by Claim \ref{c5} (a),
$r(b)<r({x_{u,v}^{k-1}}')$.
Thus $b\in X_{u,v}^{k-1}$,
which is a contradiction to $|X_{u,v}^{k-1}|=2$.
$\Box$

\bigskip

\begin{claim}\label{c8}
 The change of the interval representation of $G$ from $I$ to $I'$
creates no new bad pair $(a,b)$ such that 
$\{a,b\}\not= X_{u,v}^k$ for some $k\in [\ell_{u,v}^{\rm max}]$ and some bad pair $(u,v)$.
\end{claim}

\noindent
\textit{Proof of Claim \ref{c8}:}
For contradiction, we assume that $(a,b)$ is a new bad pair and $\{a,b\}\not= X_{u,v}^k$.
Since $(a,b)$ is a new bad pair, $I'(a)$ is a proper subset of $I(a)$.
Thus let $a\in X_{u,v}^k$ and $b\notin X_{u,v}^k$.
If $a\in X_{u,v}^k$ and $|X_{u,v}^k|=2$,
then $\ell(b)<\ell({x_{u,v}^k}')$ and $r'(a)=\ell(x_{u,v}^{k-1}) < r(b)<r({x_{u,v}^k}')$,
because of Claim \ref{c5} (a).
Thus $b\in X_{u,v}^k$, which is a contradiction.
If $a\in X_{u,v}^k$ and $|X_{u,v}^k|=1$,
then $\ell(b)<\ell({x_{u,v}^k})$ and $r'(a)=\ell(x_{u,v}^{k-1}) < r(b)<r({x_{u,v}^k})$.
Thus $b\in X_{u,v}^k$, which is the final contradiction.
$\Box$

\bigskip

In a second step, we shorten the intervals of $Y_{u,v}^{i}$ for every bad pair $(u,v)$ and $i\in [r_{u,v}^{\rm max}]$.
Let $I'':V(G)\rightarrow \mathcal{I}^{++}$ be such that $I''(y)=[r'(y_{u,v}^{k-1}),r'(y)]$ if $y\in Y_{u,v}^k$ for some bad pair $(u,v)$
and $I''(y)=I'(y)$ otherwise.
Note that bad pairs are only referred to the interval representation $I$.
Let $\ell''(x)$ and $r''(x)$ be the left and right endpoints of the interval $I''(x)$ for $x\in V(G)$, respectively.

\begin{claim}\label{c9}
$I''$ is an interval representation of $G$.
\end{claim}

\noindent
\textit{Proof of Claim \ref{c9}:}
Again, two intervals do not intersect in $I''$ if they do not intersect in $I'$ (and in $I$).
For contradiction, we assume that there are two vertices $a,b\in V(G)$ such that
$I(a)\cap I(b)\not= \emptyset$ and $I''(a)\cap I''(b)= \emptyset$.
Again, at least one interval is shortened by the change of the interval representation.
Say $a\in Y_{u,v}^k$ for some bad pair $(u,v)$ and $k\in [r_{u,v}^{\rm max}]$.

Suppose $a\in X_{w,z}^{\tilde{k}}$ for some bad pair $(w,z)$ and $\tilde{k}\in [\ell_{w,z}^{\rm max}-1]$.
By Claim \ref{c6a} (c), we have $a=x_{w,z}^{\tilde{k}}=y_{u,v}^k$.
If $y_{u,v}^{k-1}=x_{w,z}^{\tilde{k}+1}$, then we did not change the interval of $a$.
Thus we assume $y_{u,v}^{k-1}\not=x_{w,z}^{\tilde{k}+1}$.
Now $\ell(y_{u,v}^{k})<r(b)<r(y_{u,v}^{k-1})$.
The rest of the proof is similar to a symmetric version of the proof of Claim \ref{c7}.

If $a\notin X_{\tilde{u},\tilde{v}}^{\tilde{k}}$,
then $r(b)<r(y_{u,v}^{k-1})$ and $\ell({y_{u,v}^{k}}')<r(b)$, if ${y_{u,v}^{k}}'$ exists, otherwise $\ell({y_{u,v}^{k}})<r(b)$.
If $\ell(y_{u,v}^{k-1})<\ell(b)$,
then by Claim \ref{c8}, $(b,y_{u,v}^{k-1})$ is a bad pair and by Claim \ref{c5}, $I(b)=I'(b)$.
Thus Claim \ref{c1} implies the existence of a vertex, which contradicts the ``y''-version of Claim \ref{c4} (a) and (b) and hence
we suppose $\ell(b)\leq \ell(y_{u,v}^{k-1})$.
Thus $k\geq 2$, otherwise $(u',b)$ is a bad pair, which contradicts Claim \ref{c3}.
If $\ell(b)\leq \ell({y_{u,v}^{k-1}}')$, then $({y_{u,v}^{k-1}}',b)$ is a bad pair, which contradicts the ``y''-version of Claim \ref{c5} (a).
Therefore, $\ell({y_{u,v}^{k-1}}')< \ell(b)$, which implies $b\in Y_{u,v}^{k-1}$, but $b\notin \{{y_{u,v}^{k-1}},{y_{u,v}^{k-1}}'\}$, 
which contradicts the ``y''-version of Claim \ref{c4} (a).
$\Box$

\bigskip

\begin{claim}\label{c10}
The change of the interval representation of $G$ from $I$ to $I''$
creates no new bad pair $(a,b)$ such that 
$\{a,b\}\not= X_{u,v}^k$ for some $k\in [\ell_{u,v}^{\rm max}]$
or $\{a,b\}\not= Y_{u,v}^i$ for some $i\in [r_{u,v}^{\rm max}]$ and some bad pair $(u,v)$.
\end{claim}

\noindent
\textit{Proof of Claim 10:}
For contradiction, we assume that $(a,b)$ is a new bad pair and $Y_{u,v}^i\not=\{a,b\}\not= X_{u,v}^k$.
Thus $a\in X_{u,v}^k$ or $a\in Y_{u,v}^i$ and $b\notin X_{u,v}^k$ or $b\notin Y_{u,v}^i$, respectively.
If $a\in X_{u,v}^k$ and $|X_{u,v}^k|=2$,
then $\ell(b)<\ell({x_{u,v}^k}')$ and $\ell(x_{u,v}^{k-1}) < r(b)<r({x_{u,v}^k}')$.
Thus $b\in X_{u,v}^k$, which is a contradiction.
If $a\in X_{u,v}^k$ and $|X_{u,v}^k|=1$,
then $\ell(b)<\ell({x_{u,v}^k})$ and $\ell(x_{u,v}^{k-1}) < r(b)<r({x_{u,v}^k})$.
Thus $b\in X_{u,v}^k$, which is a contradiction.
If $a\in Y_{u,v}^i$ the proof is almost exactly the same.
$\Box$

\bigskip

Now we are in a position to blow up some intervals to open or half-open intervals
to get a mixed proper interval graph.
Let $I^*: V(G)\rightarrow \mathcal{I}$ be such that
\begin{align*}
	I^*(x)= \left\{
	\begin{array}{rl}
	(\ell(v),r(v)),
	& \text{if } (x,v) \text{ is a bad pair},\\
	(\ell''(x_{u,v}^k),r''(x_{u,v}^k)], 
	& \text{if } x={x_{u,v}^k}' \text{ for some bad pair } (u,v) \text{ and } k\in [\ell_{u,v}^{\rm max}-1],\\
	\left[\ell''(y_{u,v}^i),r''(y_{u,v}^i)\right), 
	& \text{if } x={y_{u,v}^i}' \text{ for some bad pair } (u,v) \text{ and } i\in [r_{u,v}^{\rm max}-1], \text{ and}\\
	\left[\ell''(x),r''(x)\right], 
	& \text{else.}
	\end{array}
	\right.
\end{align*}
Note that $I^*$ is well-defined by Claim \ref{c5} and Claim \ref{c6a}; that is,
the four cases in the definition of $I^*$ induces a partition of the vertex set of $G$.
Moreover, the interval representation $I^*$ defines a mixed proper interval graph.
As a final step, we prove that $I''$ and $I^*$ define the same graph.
Since we make every interval bigger,
we show that for every two vertices $a,b$ such that $I''(a)\cap I''(b)=\emptyset$,
we still have $I^*(a)\cap I^*(b)=\emptyset$.
For contradiction, we assume the opposite.
Let $a,b$ be two vertices such that $I''(a)\cap I''(b)=\emptyset$ and $I^*(a)\cap I^*(b)\not=\emptyset$.
It follows by our approach and definition of our interval representation $I''$,
that both $a$ and $b$ are blown up intervals.

First we suppose $a$ and $b$ are intervals that are blown up to open intervals,
that is, there are distinct vertices $\tilde{a}$ and $\tilde{b}$ such that $(a,\tilde{a})$ and $(b,\tilde{b})$ are bad pairs.
Furthermore, the intervals of $\tilde{a}$ and $\tilde{b}$ intersect not only in one point.
By Claim \ref{c2} and \ref{c3}, we assume without loss of generality, that
$\ell''(\tilde{a})<\ell''(\tilde{b})<r''(\tilde{a})<r''(\tilde{b})$.
Therefore, by the construction of $I''$,
we obtain $a$ is adjacent to $\tilde{b}$ and $\tilde{a}$ is adjacent to $b$, 
and in addition they intersect in one point, respectively.
Now, $G[\{x_{a,\tilde{a}}^1,a,\tilde{a},b,\tilde{b},y_{b,\tilde{b}}^1\}]$ is isomorphic to $T_{0,0}$, 
which is a contradiction.

Now we suppose $a$ is blown up to an open interval and $b$ is blown up to an open-closed interval
(the case closed-open is exactly symmetric).
Let $\tilde{a}$ be the vertex such that $(a,\tilde{a})$ is a bad pair.
Let $\tilde{b},u,v\in V(G)$ and $k\in \mathbb{N}$ such that
$\{b,\tilde{b}\}=X_{u,v}^k$.
We suppose $\tilde{a}\not=\tilde{b}$.
We conclude $\ell''(\tilde{a})<\ell''(\tilde{b})<r''(\tilde{a})<r''(\tilde{b})$.
As above, we conclude $a$ is adjacent to $\tilde{b}$ and $\tilde{a}$ is adjacent to $b$, 
and in addition they intersect in one point, respectively.
Thus $G[\{x_{a,\tilde{a}}^1,a,\tilde{a},v,u,y_{u,v}^1\}\cup \bigcup_{i=1}^k X_{u,v}^i]$ induces a $T_{k,0}$,
which is a contradiction.
Now we suppose $\tilde{a}=\tilde{b}$.
We conclude that $G[\{x_{a,\tilde{a}}^1,a,v,u,y_{u,v}^1\}\cup \bigcup_{i=1}^{k} X_{u,v}^i]$
is isomorphic to $R_k$,
which is a contradiction.

It is easy to see that $a$ and $b$ cannot be both blown up to closed-open or both open-closed intervals,
because $G$ is $R_k$-free for $k\geq 0$ and the definition of $I''$.

Therefore, we consider finally the case that $a$ is blown up to a closed-open and $b$ to an open-closed interval.
Let $\tilde{a},\tilde{b},u,v,w,z \in V(G)$ and $k,\tilde{k}\in \mathbb{N}$
such that $\{a,\tilde{a}\}=Y_{u,v}^k$ and $\{b,\tilde{b}\}=X_{w,z}^{\tilde{k}}$.
First we suppose $\tilde{a}\not=\tilde{b}$.
Again, we obtain $\ell''(\tilde{a})<\ell''(\tilde{b})<r''(\tilde{a})<r''(\tilde{b})$
and $a$ is adjacent to $\tilde{b}$ and $\tilde{a}$ is adjacent to $b$, 
and furthermore they intersect in one point, respectively.
Thus $G[\{x_{u,v}^1,u,v,w,z,y_{w,z}^1\}\cup \bigcup_{i=1}^k Y_{u,v}^i\cup \bigcup_{i=1}^{\tilde{k}} X_{w,z}^i]$ is isomorphic to $T_{k,\tilde{k}}$.
Next we suppose $\tilde{a}=\tilde{b}$ and hence 
$G[\{x_{u,v}^1,u,v,w,z,y_{w,z}^1\}\cup \bigcup_{i=1}^k Y_{u,v}^i\cup \bigcup_{i=1}^{\tilde{k}} X_{w,z}^i]$ is isomorphic to $R_{k+\tilde{k}}$.
This is the final contradiction and completes the proof of Theorem \ref{mainthm}.
$\Box$

\bigskip

\begin{figure}[t]
\begin{center}
\begin{tikzpicture}[scale=1]
\def\ver{0.1} 
\def\x{1}

\def\xa{10}
\def\ya{0}

\def\xb{0}
\def\yb{0}

\def\xc{0}
\def\yc{0}

\def\xd{5}
\def\yd{0}

\path[fill] (\xc,\yc) circle (\ver);
\path[fill] (\xc+1,\yc) circle (\ver);
\path[fill] (\xc+2,\yc) circle (\ver);
\path[fill] (\xc+3,\yc) circle (\ver);
\path[fill] (\xc+0.5,\yc+1) circle (\ver);
\path[fill] (\xc+1.5,\yc+1) circle (\ver);
\path[fill] (\xc+2.5,\yc+1) circle (\ver);
\path[fill] (\xc+1,\yc+1.5) circle (\ver);

\draw[thick] (\xc,\yc)--(\xc+1,\yc)--(\xc+2,\yc)--(\xc+3,\yc)
(\xc+0.5,\yc+1)--(\xc+1,\yc)--(\xc+1.5,\yc+1)--(\xc+2,\yc)--(\xc+2.5,\yc+1)
(\xc,\yc)--(\xc+0.5,\yc+1)--(\xc+1,\yc+1.5)--(\xc+1.5,\yc+1)
(\xc+1,\yc)--(\xc+1,\yc+1.5);

\draw[thick] (\xc+2,\yc) arc(-20:85:1.16);

\node (1) at (\xc+1.5,\yc-0.7) {$S_2''$};

\path[fill] (\xd-1,\yd) circle (\ver);
\path[fill] (\xd-0.5,\yd+1) circle (\ver);
\path[fill] (\xd,\yd) circle (\ver);
\path[fill] (\xd+1,\yd) circle (\ver);
\path[fill] (\xd+2,\yd) circle (\ver);
\path[fill] (\xd+3,\yd) circle (\ver);
\path[fill] (\xd+0.5,\yd+1) circle (\ver);
\path[fill] (\xd+1.5,\yd+1) circle (\ver);
\path[fill] (\xd+2.5,\yd+1) circle (\ver);
\path[fill] (\xd,\yd+1.5) circle (\ver);

\draw[thick] (\xd-1,\yd)--(\xd+3,\yd)
(\xd-1,\yd)--(\xd-0.5,\yd+1)--(\xd,\yd)--(\xd+0.5,\yd+1)--(\xd+1,\yd)--(\xd+1.5,\yd+1)--(\xd+2,\yd)--(\xd+2.5,\yd+1)
(\xd-0.5,\yd+1)--(\xd,\yd+1.5)--(\xd+0.5,\yd+1)
(\xd,\yd)--(\xd,\yd+1.5);

\draw[thick] (\xd+1,\yd) arc(-20:85:1.16);

\node (1) at (\xd+1,\yd-0.7) {$S_3''$};

\path[fill] (\xa-1,\ya) circle (\ver);
\path[fill] (\xa,\ya) circle (\ver);
\path[fill] (\xa+1,\ya) circle (\ver);
\path[fill] (\xa+2,\ya) circle (\ver);
\path[fill] (\xa+3.5,\ya) circle (\ver);
\path[fill] (\xa+4.5,\ya) circle (\ver);
\path[fill] (\xa+5.5,\ya) circle (\ver);
\path[fill] (\xa-0.5,\ya+1) circle (\ver);
\path[fill] (\xa+0.5,\ya+1) circle (\ver);
\path[fill] (\xa+1.5,\ya+1) circle (\ver);
\path[fill] (\xa+4,\ya+1) circle (\ver);
\path[fill] (\xa+5,\ya+1) circle (\ver);
\path[fill] (\xa,\ya+1.5) circle (\ver);

\draw[thick] (\xa-1,\ya)--(\xa+2,\ya)
(\xa+3.5,\ya)--(\xa+5.5,\ya)
(\xa-1,\ya)--(\xa-0.5,\ya+1)--(\xa,\ya)--(\xa+0.5,\ya+1)--(\xa+1,\ya)--(\xa+1.5,\ya+1)--(\xa+2,\ya)
(\xa+3.5,\ya)--(\xa+4,\ya+1)--(\xa+4.5,\ya)--(\xa+5,\ya+1)
(\xa-0.5,\ya+1)--(\xa,\ya+1.5)--(\xa+0.5,\ya+1)
(\xa,\ya)--(\xa,\ya+1.5);

\draw[thick] (\xa+1,\ya) arc(-20:85:1.16);

\node (1) at (\xa+5,\ya-0.7) {$S_i''$};

\fill (\xa+2.35,\ya) circle (\ver/2);
\fill (\xa+2.75,\ya) circle (\ver/2);
\fill (\xa+3.15,\ya) circle (\ver/2);

\draw[thick,decoration={brace,mirror,raise=0.2cm},decorate] (\xa-1,\ya) -- (\xa+4.5,\ya)
node [pos=0.5,anchor=north,yshift=-0.4cm] {$i$ triangles};

\end{tikzpicture}
\end{center}
\caption{The class $\mathcal{S}_i''$.}\label{graphsS''}
\end{figure}

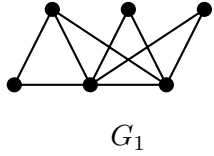
\begin{figure}[t]
\begin{center}
\begin{tikzpicture}[scale=1]
\def\ver{0.1} 
\def\x{1}

\def\xa{0}
\def\ya{0}

\path[fill] (\xa+1,\ya) circle (\ver);
\path[fill] (\xa+2,\ya) circle (\ver);
\path[fill] (\xa+0.5,\ya+1) circle (\ver);
\path[fill] (\xa+1.5,\ya+1) circle (\ver);
\path[fill] (\xa+2.5,\ya+1) circle (\ver);
\fill (\xa,\ya) circle (\ver);

\draw[thick] (\xa+1,\ya)--(\xa+2,\ya)
(\xa+1,\ya)--(\xa+0.5,\ya+1)--(\xa+2,\ya)
(\xa+1,\ya)--(\xa+2.5,\ya+1)--(\xa+2,\ya)
(\xa+1,\ya)--(\xa+1.5,\ya+1)--(\xa+2,\ya)
(\xa+1,\ya)--(\xa,\ya)--(\xa+0.5,\ya+1);

\node (1) at (\xa+1.5,\ya-0.7) {$G_1$};

\end{tikzpicture}
\end{center}
\caption{The graph $G_1$.}\label{graphG1}
\end{figure}

In Theorem \ref{mainthm} we only consider twin-free $\mathcal{U}$-graphs to reduce the number of case distinctions in the proof.
In Corollary \ref{maincoro} we resolve this technical condition.
See Figure \ref{graphsS''} and \ref{graphG1} for illustration.
Let $\mathcal{S''}=\bigcup_{i=2}^{\infty}\{S_i''\}$.

\begin{coro}\label{maincoro}
A graph $G$ is a mixed unit interval graph if and only if 
$G$ is a $\{G_1\}\cup\mathcal{R}\cup \mathcal{S}\cup \mathcal{S''}\cup \mathcal{T}$-free interval graph.
\end{coro}

\noindent
\textit{Proof of Corollary \ref{maincoro}}:
We first show that $ \{G_1\}\cup\mathcal{R}\cup \mathcal{S}\cup \mathcal{S''}\cup \mathcal{T}$
is the set of all twin-free graphs that contain all graphs of
$\{K_{2,3}^*\}\cup\mathcal{R}\cup \mathcal{S}\cup \mathcal{S'}\cup \mathcal{T}$
and are minimal with subject to induced subgraphs.
We leave it as an exercise to show that $G_1$ is the only minimal twin-free and $R_0$-free graph that contains $K_{2,3}^*$.
Since all graphs in $\mathcal{R}\cup \mathcal{S}\cup \mathcal{T}$ are twin-free graphs,
there is nothing to show.

Let now $G\in \mathcal{S'}$, that is $G=S_k'$ for some $k\in \mathbb{N}$.
With the notation as in the proof of Theorem \ref{mainthm},
$G$ can be interpreted as a bad pair $(u,v)$ together with $\{y_{u,v}^1\}\cup\bigcup_{i=1}^k X_{u,v}^i$
such that $|X_{u,v}^i|=2$ if $i<k$ and $|X_{u,v}^k|=3$.
Note that Claim \ref{c4} (b) of Theorem \ref{mainthm} is still true even if $G$ is not $\mathcal{S}'$-free.
Therefore, we know that the vertices in $X_{u,v}^i$ cannot be distinguished by vertices from the right.
Thus the vertices that distinguish the vertices in $X_{u,v}^k$ are only adjacent to $X_{u,v}^k$.
Clearly, there are at least two of them, say $a,b$.
Without loss of generality $a$ and $b$ they do not have the same neighborhood on $X_{u,v}^k$.
We conclude either $N_{G[X_{u,v}^k]}(a)\subset N_{G[X_{u,v}^k]}(b)$ or $N_{G[X_{u,v}^k]}(b)\subset N_{G[X_{u,v}^k]}(a)$.
We assume the first possibility.
Since $0<|N_{G[X_{u,v}^k]}(x)\cap X_{u,v}^k|<3$ for $x\in \{a,b\}$,
it follows $|N_{G[X_{u,v}^k]}(a)\cap X_{u,v}^k|=1$ and $|N_{G[X_{u,v}^k]}(b)\cap X_{u,v}^k|=2$.
Since $G$ is $R_k$-free, $a$ and $b$ are adjacent.
Now $G[\bigcup_{i=1}^k X_{u,v}^i\cup \{a,b,u,v,y_{u,v}^1\}]$ is isomorphic to $S_{k+1}''$.
This completes this part of the proof.

Let $G$ be an interval graph.
The relation $\sim$, where $u\sim v$ if and only if $u$ and $v$ are twins, defines an equivalence relation on $V(G)$.
Let $U\subseteq V(G)$ such that there is exactly one vertex of every equivalence class in $U$.
Therefore, $G[U]$ is a twin-free graph.
Furthermore, $G$ contains an induced subgraph in $\{G_1\}\cup\mathcal{R}\cup \mathcal{S}\cup \mathcal{S''}\cup \mathcal{T}$
if and only if $G[U]$ contains an induced subgraph in $\{K_{2,3}^*\}\cup\mathcal{R}\cup \mathcal{S}\cup \mathcal{S'}\cup \mathcal{T}$.
In addition, $G[U]$ is a twin-free $\mathcal{U}$-graph if and only if $G$ is a $\mathcal{U}$-graph.
By Theorem \ref{mainthm} this completes the proof.
$\Box$


\begin{thebibliography}{1}

\bibitem{corneil}
D.G. Corneil, A simple 3-sweep LBFS algorithm for the recognition of unit interval graphs, Discrete Appl. Math. \textbf{138}, 371-379 (2004).

\bibitem{corneiletal}
D.G. Corneil, S. Olariu and L. Stewart, The ultimate interval graph recognition algorithm?,
Proceedings of the {N}inth {A}nnual {ACM}-{SIAM} {S}ymposium on {D}iscrete {A}lgorithms ({S}an {F}rancisco, {CA}, 1998), 175-180.

\bibitem{dlprs}
M.C. Dourado, V.B. Le, F. Protti, D. Rautenbach and J.L. Szwarcfiter, Mixed unit interval graphs, Discrete Math. \textbf{312}, 3357-3363 (2012).

\bibitem{fishburn}
P.C. Fishburn. Interval Orders and Interval Graphs. John Wiley \& Sons (1985).

\bibitem{fm}
P. Frankl and H. Maehara, Open interval-graphs versus closed interval-graphs, Discrete Math. \textbf{63}, 97-100 (1987).

\bibitem{golumbic}
M.C. Golumbic, Algorithmic Graph Theory and Perfect Graphs, vol. 57. Annals of Discrete Mathematics, Amsterdam, The Netherlands, 2004.

\bibitem{lr}
V.B. Le and D. Rautenbach, Integral Mixed Unit Interval Graphs, Lecture Notes in Computer Science 7434, 495-506 (2012).

\bibitem{LekBo}
C.G. Lekkerkerker and J.C. Boland, Representation of a finite graph by a set of intervals on the real line, Fund. Math., \textbf{51}, 45-64 (1962).

\bibitem{rs}
D. Rautenbach and J.L. Szwarcfiter, Unit Interval Graphs of Open and Closed Intervals, J. Graph Theory \textbf{72}(4), 418-429 (2013).

\bibitem{roberts}
F.S. Roberts, Indifference graphs, in F. Harary (Ed.), Proof Techniques in Graph Theory, Academic Press, pp.139-146 (1969).

\end{thebibliography}
\end{document}